\documentclass[11pt]{amsart}

\usepackage[foot]{amsaddr}

\usepackage[T1]{fontenc}
\usepackage[latin9]{inputenc}
\usepackage{enumerate}
\usepackage{endnotes}
\usepackage{amsthm}
\usepackage{amstext}
\usepackage{graphicx}
\usepackage{amssymb}
\usepackage{graphicx}
\usepackage{standalone}
\usepackage[stable]{footmisc}
\usepackage{lmodern}
\usepackage{fullpage}
\usepackage{xcolor}
\usepackage{comment}
\graphicspath{{./Images/}}
\usepackage{amsmath} % for substack
\usepackage{amssymb}
\usepackage[T1]{fontenc}    % use 8-bit T1 fonts
\usepackage{hyperref}       % hyperlinks
\hypersetup{final}
\usepackage{pgfplots}
\pgfplotsset{compat=1.16}
\usepackage{url}            % simple URL typesetting
\usepackage{float}
\usepackage{multirow}
\usepackage{booktabs}       % professional-quality tables
\usepackage{amsfonts}       % blackboard math symbols
\usepackage{nicefrac}       % compact symbols for 1/2, etc.
\usepackage{lipsum}
\usepackage{verbatim}
\usepackage{tikz}
\usepackage{caption}
\usepackage{subcaption}
\usepackage{blindtext}
\usepackage{endnotes}
\usepackage{colortbl}
\usepackage{algorithm}
\usepackage{algpseudocode}
\usepackage{subfiles} % Best loaded last in the preamble
\newtheorem{theorem}{Theorem}[section]
\newtheorem{definition}{Definition}[section]
\newtheorem{lemma}[theorem]{Lemma}
\theoremstyle{remark}
\newtheorem{remark}[theorem]{Remark}

\def\keywordname{{\bfseries \emph{Keywords}}}%
\def\keywords#1{\par\addvspace\medskipamount{\rightskip=0pt plus1cm
		\def\and{\ifhmode\unskip\nobreak\fi\ $\cdot$
		}\noindent\keywordname\enspace\ignorespaces#1\par}}
	
\newcommand{\Rb}{\mathbb{R}}
\newcommand{\Cb}{\mathbb{C}}
\newcommand{\Zb}{\mathbb{Z}}

\newcommand{\rect}[3]{\draw [ultra thick,fill=#1] (#2) rectangle (#3);}

\newcommand{\ceil}[1]{
	\left\lceil #1 \right \rceil
}

\newcommand{\abs}[1]{
	\left\lvert #1 \right\rvert
}

\newcommand{\magn}[1]{
	\left\lVert #1 \right\rVert
}

\newcommand{\bkt}[1]{
	\left(#1\right)
}

\newcommand{\dsum}{\displaystyle \sum}

\newcommand{\dcap}{\displaystyle \cap}
\newcommand{\dbcup}{\displaystyle \bigcup}

\numberwithin{equation}{section}
\begin{document}
	
	\title{HODLR2D: A new class of hierarchical matrices}
	
	\author{V A Kandappan}
	\email{kandappanva@gmail.com}
	
	\author{Vaishnavi Gujjula}
%	\address{Department of Mathematics, Indian Institute of Technology Madras}
	\email{vaishnavihp@gmail.com}
	
	\author{Sivaram Ambikasaran}

	\email{sivaambi@smail.iitm.ac.in}
	\address{Department of Mathematics, Indian Institute of Technology Madras}

\maketitle

\begin{abstract}
    This article introduces HODLR2D, a new hierarchical low-rank representation for a class of dense matrices arising out of $N$ body problems in two dimensions. Using this new hierarchical framework, we propose a new fast matrix-vector product that scales almost linearly. We apply this fast matrix-vector product to accelerate the iterative solution of large dense linear systems arising out of radial basis function interpolation and discretized integral equation. The space and computational complexity of HODLR2D matrix-vector products scales as $\mathcal{O}(pN \log(N))$, where $p$ is the maximum rank of the compressed matrix subblocks. We prove that $p \in \mathcal{O}\bkt{\log\bkt{N}\log\bkt{\log\bkt{N}}}$, which ensures that the storage and computational complexity of HODLR2D matrix-vector products remain tractable for large $N$. Additionally, we also present the parallel scalability of HODLR2D as part of this article.\\
    \smallskip 
    \keywords{\textbf{Key words.} Hierarchical matrices, Low-rank approximations, N-body problems, Iterative methods, Radial basis functions} 
\end{abstract}

	\section{Introduction}

This article considers a class of dense rank structured matrices that have a hierarchical low-rank structure~\cite{hmat_1999,hmat_2003}. Hierarchical matrices frequently arise in applications such as radial basis function interpolation~\cite{gumerov2007fast},  electromagnetic scattering~\cite{ROKHLIN1990414}, geostatistics~\cite{ambikasaran2013large}, machine learning~\cite{gray2000n}, Gaussian process regression~\cite{ambikasaran2016fast}, etc. In this article, we focus on a new class of hierarchical matrices arising out of two-dimensional problems. We term this new class of Hierarchical matrices as HODLR2D. For this new class of hierarchical matrices, HODLR2D, we propose a fast matrix-vector product that scales as $\mathcal{O}\bkt{p N \log N}$, where $p$ is the maximum rank of the compressed matrix submatrices. Further, we prove that the rank $p$ can only grow at most as $\mathcal{O}\bkt{\log N \log \bkt{\log N}}$.

Fast matrix-vector products for $N$ body problems have been studied even before the advent of hierarchical matrices. Chronologically, Barnes and Hut~\cite{treecode} reduced the computational cost of the matrix-vector product arising out of three dimensional N-body problem from $\mathcal{O}(N^2)$ to $\mathcal{O}(N\log (N))$ by using an oct tree to hierarchically sub-divide the domain and then efficiently approximating "far-away" interactions. The Fast Multipole Method (FMM)~\cite{fmm_1987} brought the cost for matrix-vector product further down to $\mathcal{O}(N)$. The literature on the Barnes and Hut algorithm~\cite{treecode} and FMM~\cite{fmm_1987} is extensive, and we direct readers to articles~\cite{panel,Darve_2000, Fong_2009,FMMcourse,FMMtypes,martinsson2007accelerated} and the references therein. Although these methods are specific to the N-body problem, they form the foundation for the hierarchical matrices (from now on, termed as $\mathcal{H}$-matrices). Hackbusch et al.~\cite{hmat_1999,hmat_2003} formally defined these $\mathcal{H}$-matrices by hierarchically sub-dividing the underlying matrix using an appropriate tree structure and representing certain off-diagonal blocks at different levels in the tree, based on certain admissibility criterion, as low-rank matrices. The Hierarchically Off-Diagonal Low Rank (HODLR)~\cite{SA_FDS_2013} matrix relies on hierarchically sub-dividing the matrix using a binary tree and represents all off-diagonal blocks at each level in the tree as low-rank matrices. The computational cost of matrix-vector product using HODLR matrix scales as $\mathcal{O}(pN\log (N))$, where $p$ is the maximum rank of the off-diagonal blocks. One disadvantage of HODLR matrices is that $p$ remains almost constant only for $1$D problems and grows significantly with $N$ for higher dimensions (roughly as $\mathcal{O}\bkt{\sqrt{N}}$ in $2$D and as $\mathcal{O}\bkt{\sqrt[3]{N^2}}$ in $3$D)~\cite{ambikasaran2014inverse}. There are also other hierarchical matrix structures such as HSS~\cite{chandrasekaran2002fast}, HBS~\cite{gillman2012direct}, $\mathcal{H}^2$~\cite{h2mat}, etc. These different hierarchical matrix structures differ based on the following three criteria.
\begin{enumerate}
	\item Choice of tree structure to sub-divide the matrix
	\item Identifying submatrices that are efficiently represented as low-rank matrices
	\item Whether the row and column basis of the low-rank submatrices are nested or not (i.e., whether the row and column basis for the identified low-rank submatrices can be constructed from the row and column basis of its children)~\cite{ncaBeb,ncaZhao,ncaVaish}
\end{enumerate}

The new HODLR2D matrix representation for $N$-body problems in $2$D is based on the following choices.
\begin{enumerate}
	\item A quadtree is used to sub-divide the underlying two dimensional domain, which in turn sub-divides the matrix
	\item The matrix corresponding to the interaction between the sub domains that share a vertex is efficiently represented as a low-rank matrix
	\item The row and column basis of the low-rank submatrices are \textbf{not nested}
	\item The low-rank representation of the desired low-rank submatrices are obtained using Adaptive Cross Approximation (from now on abbreviated as ACA)~\cite{aca_2003,aca_zhao}.
\end{enumerate}
These are discussed in detail in Sections~\ref{sec4}. Once the HODLR2D representation is obtained the computational cost to perform matrix-vector product scales as $\mathcal{O}\bkt{pN\log N}$, where $p$ is the rank of compressed submatrices and scales $\mathcal{O}\bkt{\log N \log \log N}$. A detailed analysis, including the proof that $p \in \mathcal{O}\bkt{\log N \log \log N}$, is presented in Section~\ref{sec3}.
Section~\ref{sec5} uses this fast HODLR2D matrix-vector product to accelerate iterative solver for dense linear systems arising out of radial basis function interpolation and a discretized integral equation. Finally, Section~\ref{sec6} illustrates the parallel scalability of these HODLR2D matrix-vector products.

The main highlights of this article are as follows:
\begin{itemize}
	\item New class of Hierarchical matrix (HODLR2D) for $N$-body problems in $2$D is proposed.
	\item Bounds on rank of interaction between neighboring sub-domains in $2$D is proved.
	\item Storage and computational complexity of HODLR2D matrix-vector products scales as $\mathcal{O}\bkt{pN\log\bkt{N}}$, where $p\in \mathcal{O}\bkt{\log\bkt{N}\log\bkt{\log\bkt{N}}}$.
	\item Fast HODLR2D matrix-vector product is leveraged to accelerate dense matrix solvers and compared against HODLR and $\mathcal{H}$-matrix with standard admissibility criterion.
	\item HODLR2D is not only significantly better than HODLR but also provides an attractive alternative to $\mathcal{H}$-matrices with standard admissibility criterion for $2$D problems.
	\item Parallel scalability of HODLR2D matrix-vector product is studied.
\end{itemize}
	\section{Preliminaries}
	\label{sec2}
	We begin by establishing some notations using two hierarchical low-rank representations: (i) A generic $\mathcal{H}$ matrix; (ii) HODLR matrix. We look at the above two structures since we will be comparing our HODLR2D against these two hierarchical matrices.
	
	\subsection{Notations}
	Let $B \subset \Rb^2$ be a box containing $N$ particles whose interaction is given by the matrix $A \in \Rb^{n \times n}$, i.e., $A_{ij}$ is the interaction between the $i^{th}$ and $j^{th}$ particle (We will assume that the matrix is symmetric for pedagogical reasons). Let $\mathcal{T}^{L}$ represent a $L$ level tree, which subdivides the box $B$ hierarchically. Each node in the tree represents a set of particles inside the box $B$, and $\mathcal{I}$ be the index set that maps all these particles in the box B. The following definitions help define $\mathcal{T}^{L}$.
	
	\begin{definition}
		$\mathcal{N}_{i}^{(l)}$ denote the $i^{th}$ node at level $l$ in the tree $\mathcal{T}^{L}$.
	\end{definition}
	\begin{definition}
		$\mathcal{I}_{i}^{(l)}$ is the \textbf{index set} of particles corresponding to $\mathcal{N}_{i}^{(l)}$.
	\end{definition}
	\begin{definition}
		The cluster $\mathcal{C}_i^{(l)}$ contains all the coordinate pairs corresponding to $\mathcal{N}_{i}^{(l)}$.
	\end{definition}
	For example, in the box $B$ considered $\mathcal{N}_0^{(0)}$ contains $\mathcal{I}_0^{(0)}$ and $\mathcal{C}_0^{(0)}$, where $\mathcal{I}_0^{(0)} = \{1,2,\ldots,N\}$ and $\mathcal{C}_0^{(0)} = \{(x_1,y_1),(x_2,y_2),(x_3,y_3),\ldots,(x_N,y_N)\}$, where $\bkt{x_i,y_i}$ is the location of the $i^{th}$ particle.
	
	$\mathcal{N}_0^{(0)}$ is called the root of $\mathcal{T}^{L}$.
	\begin{definition}
		Consider two nodes $\mathcal{N}^{(l)}_{i}$ at level $l$ and $\mathcal{N}^{(l+1)}_{j}$ at level $l+1$. If $\mathcal{I}^{(l+1)}_{j} \subseteq \mathcal{I}^{(l)}_{i}$, then $\mathcal{N}^{(l)}_{i}$ is termed as the parent of $\mathcal{N}^{(l+1)}_{j}$  in $\mathcal{T}^{L}$ and $\mathcal{N}^{(l+1)}_{j}$ is the child of $\mathcal{N}^{(l)}_{i}$ in $\mathcal{T}^{L}$. Note that we immediately have $\mathcal{C}^{(l+1)}_{j} \subseteq \mathcal{C}^{(l)}_{i}$. Further, the index set $\mathcal{I}^{(l)}_{i} = \bigcup_{\mathcal{N}^{(l+1)}_{k} \in child(\mathcal{N}^{(l)}_{i})} \mathcal{I}^{(l+1)}_{k}$
	\end{definition}
	\begin{definition}
		Consider two nodes $\mathcal{N}^{(l)}_{i}$ and $\mathcal{N}^{(l)}_{j}$ at level $l$ in $\mathcal{T}^{L}$. If $\mathcal{N}^{(l)}_{i}$  and $\mathcal{N}^{(l)}_{j}$ have same parent, then $\mathcal{N}^{(l)}_{i}$ is the sibling of $\mathcal{N}^{(l)}_{j}$ in $\mathcal{T}^{L}$.
	\end{definition}
	\begin{definition}
		For a node $\mathcal{N}$, if child($\mathcal{N}$) = $\emptyset$, then $\mathcal{N}$ is the leaf in $\mathcal{T}^{L}$.
	\end{definition}
	The interaction between the points in the same cluster with index set $I$ is denoted by the matrix block $A(I,I)$; the interaction of particles in the cluster with index set $J$ with the particles in the cluster with index set $I$ is $A(I,J)$. Throughout the article, we use the following definition to compute the numerical rank of a matrix.
	\begin{definition}
		\textbf{Numerical Rank of a matrix} Given $\epsilon>0$, the $\epsilon$-rank of the matrix $A \in \Cb^{n \times n}$, denoted by $r_{\epsilon}(A)$, is given by
		$$r_{\epsilon}(A) = \max \left\{k \in \{1,2,\ldots,N\}: \dfrac{\sigma_{k}}{\sigma_1} > \epsilon \right\}$$ where $\sigma_1\geq \sigma_2 \geq ... \geq \sigma_N \geq 0$ are the singular value of the matrix $A$.
	\end{definition}
	\begin{definition}
		Let $A \in \Cb^{m \times n}$. We say that a matrix algorithm on the matrix $A$ scales \textbf{almost linearly} if the computational complexity of the matrix algorithm (measured in terms of flop counts) scales as $\mathcal{O}\bkt{\bkt{m+n}^{1+\epsilon}}$ for all $\epsilon > 0$.
	\end{definition}
	
	\subsection{\texorpdfstring{$\mathcal{H}$}{Hc}-matrix}
	The hierarchical low-rank matrices ($\mathcal{H}$-matrix) operate on the tree $\mathcal{T}^L$ (like $2^d$-tree for $d$ dimensions, a K-D tree, etc.) to represent the clusters formed by hierarchical subdivision of the underlying domain. The submatrix corresponding to the interaction between two nodes at the same level in the hierarchical tree can be efficiently approximated using a low-rank matrix, if it agrees with the admissibility condition as given below. Let us consider two nodes at level $l$, $\mathcal{N}_{i}^{(l)}$ and $\mathcal{N}_{j}^{(l)}$. The admissibility condition for $\mathcal{N}_{j}^{(l)}$ to $\mathcal{N}_{i}^{(l)}$ is given by \eqref{eq:adm}.
	\begin{equation}\label{eq:adm}
		min\bkt{diam\bkt{\mathcal{C}_{i}^{(l)}},diam\bkt{\mathcal{C}_{j}^{(l)}}} \leq \eta dist\bkt{\mathcal{C}_{i}^{(l)},\mathcal{C}_{j}^{(l)}}
	\end{equation}
	where $diam(\mathcal{C})$ is the Euclidean diameter of the cluster $\mathcal{C}$ and $dist\bkt{\mathcal{C}_{i}^{(l)},\mathcal{C}_{j}^{(l)}}$ is the Euclidean distance between the two clusters. If the clusters $\mathcal{C}_{i}^{(l)}$ and $\mathcal{C}_{j}^{(l)}$ satisfy \eqref{eq:adm}, then they are defined as admissible clusters. We direct our readers to \cite{hmat_2003} for a detailed description of $\mathcal{H}$-matrix and its practical implementation details.
	
	\textbf{\emph{In this article, when we consider $\mathcal{H}$-matrix, we use a quadtree to subdivide the domain and the admissibility condition we use sets the value of $\eta$ to be $\sqrt{2}$. This $\mathcal{H}$-matrix representation is one of the baseline representation that we use to compare the performance of our proposed HODLR2D representation}}.
	
	\subsection{HODLR matrix}
	We now describe the HODLR~\cite{ambikasaran2013large} matrix. The HODLR matrix approximates the submatrix of interactions between any two disjoint clusters as a low-rank matrix. Typically, the HODLR matrix subdivides the domain using a K-D tree. Following the conventions in~\cite{SA_FDS_2013}, the 1-level HODLR representation of the dense matrix $A \in \mathbb{R}^{N \times N}$ is given in Equation~\eqref{eq:hodlr}.
	\begin{equation}
		A = K_{1}^{(0)} = \begin{bmatrix}
			K_{1}^{(1)} &
			K_{12}^{(1)}\\
			K_{21}^{(1)} &
			K_{2}^{(1)} 
		\end{bmatrix} = \begin{bmatrix}
			K_{1}^{(1)} &
			U^{(1)}_{1} V^{(1)^T}_{2}\\
			U^{(1)}_{2} V^{(1)^T}_{1} &
			K_{2}^{(1)} 
		\end{bmatrix}\label{eq:hodlr}
	\end{equation}
	where $K_1^{(1)}, K_2^{(1)} \in \Rb^{N/2 \times N/2}$, $U_1,U_2,V_1,V_2 \in \Rb^{N/2 \times p}$ and $p \ll N$.
	
\begin{figure}[H]
		\begin{subfigure}[b]{0.3\textwidth}
			\centering
			\begin{tikzpicture}
				[%%%%%%%%%%%%%%%%%%%%%%%%%%%%%%
				box/.style={rectangle,draw=black, minimum size=0.5cm},
				]%%%%%%%%%%%%%%%%%%%%%%%%%%%%%%
				
				\draw[fill=cyan!50] (0,0) rectangle +(2,2);
				\draw[fill=cyan!50] (2,2) rectangle +(2,2);
				\draw[fill=red] (0,2) rectangle +(2,2);
				\draw[fill=red] (2,0) rectangle +(2,2);
			\end{tikzpicture}
			\caption{level = 1}
		\end{subfigure}%
		\hfill
		\begin{subfigure}[b]{0.3\textwidth}
			\centering
			\begin{tikzpicture}
				[%%%%%%%%%%%%%%%%%%%%%%%%%%%%%%
				box/.style={rectangle,draw=black, minimum size=0.5cm},
				]%%%%%%%%%%%%%%%%%%%%%%%%%%%%%%
				
				\draw[fill=cyan!50] (0,0) rectangle +(2,2);
				\draw[fill=cyan!50] (2,2) rectangle +(2,2);
				\draw[fill=cyan!50] (0,2) rectangle +(1,1);
				\draw[fill=cyan!50] (1,3) rectangle +(1,1);
				\draw[fill=cyan!50] (3,1) rectangle +(1,1);
				\draw[fill=red] (1,2) rectangle +(1,1);
				\draw[fill=red] (0,3) rectangle +(1,1);
				\draw[fill=red] (3,0) rectangle +(1,1);
				\draw[fill=cyan!50] (2,0) rectangle +(1,1);
				\draw[fill=red] (2,1) rectangle +(1,1);
			\end{tikzpicture}
			\caption{level = 2}
		\end{subfigure}%
		\hfill
		\begin{subfigure}[b]{0.3\textwidth}
			\centering
			\begin{tikzpicture}
				[%%%%%%%%%%%%%%%%%%%%%%%%%%%%%%
				box/.style={rectangle,draw=black, minimum size=0.5cm},
				]%%%%%%%%%%%%%%%%%%%%%%%%%%%%%%
				\draw[fill=cyan!50] (0,0) rectangle +(2,2);
				\draw[fill=cyan!50] (0,2) rectangle +(1,1);
				\draw[fill=cyan!50] (0,3) rectangle +(0.5,0.5);
				\draw[fill=red] (0,3.5) rectangle (0.5,4);
				\draw[fill=cyan!50] (0.5,3.5) rectangle (1,4);
				\draw[fill=red] (0.5,3) rectangle (1,3.5);
				\draw[fill=cyan!50] (1,3) rectangle (2,4);
				\draw[fill=cyan!50] (1,2) rectangle +(0.5,0.5);
				\draw[fill=red] (1,2.5) rectangle (1.5,3);
				\draw[fill=cyan!50] (1.5,2.5) rectangle (2,3);
				\draw[fill=red] (1.5,2) rectangle (2,2.5);
				\draw[fill=cyan!50] (2,2) rectangle (4,4);
				\draw[fill=cyan!50] (2,0) rectangle +(1,1);
				\draw[fill=cyan!50] (2,1) rectangle +(0.5,0.5);
				\draw[fill=red] (2,1.5) rectangle (2.5,2);
				\draw[fill=cyan!50] (2.5,1.5) rectangle (3,2);
				\draw[fill=red] (2.5,1) rectangle (3,1.5);
				\draw[fill=cyan!50] (3,1) rectangle (4,2);
				\draw[fill=cyan!50] (3,0) rectangle +(0.5,0.5);
				\draw[fill=red] (3,0.5) rectangle (3.5,1);
				\draw[fill=cyan!50] (3.5,0.5) rectangle (4,1);
				\draw[fill=red] (3.5,0) rectangle (4,0.5);
			\end{tikzpicture}
			\caption{level = 3}
		\end{subfigure}%
		\hfill
		\begin{subfigure}[b]{0.4\textwidth}
			\centering
			\begin{tikzpicture}
				[%%%%%%%%%%%%%%%%%%%%%%%%%%%%%%
				box/.style={rectangle,draw=black, minimum size=0.5cm},
				]%%%%%%%%%%%%%%%%%%%%%%%%%%%%%%
				
				\node[box,fill=red,label=right:Full-rank Matrix (Self Interaction),anchor=west] at (5,3){};
				\node[box,fill=cyan,label=right:Low-rank representation ($UV^T$),anchor=west] at (5,2){};
			\end{tikzpicture}
			\caption*{}
		\end{subfigure}%
		\caption{A HODLR matrix at different levels.}
		\label{fig:hodlr_mat}
	\end{figure}
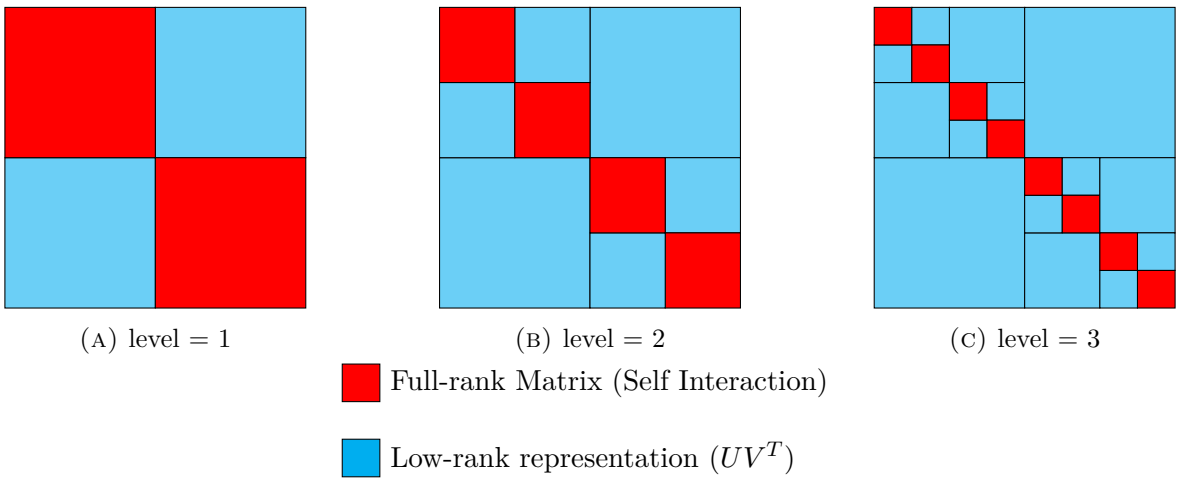

	The dense block matrices $K_1^{(1)}$, $K_2^{(1)}$ can further be represented as a HODLR matrix. In general, for an $L$ level HODLR matrix, its $i^{th}$ diagonal block at a level $l$, $K_{i}^{(l)}$ is given by Equation~\eqref{eq:hodlr_gen} where $1 \leq i \leq 2^l$  and $0\leq l \leq L$.
	\begin{equation}
		K_{i}^{(l)} = \begin{bmatrix}
			K_{2i-1}^{(l+1)} &
			U^{(l+1)}_{2i-1} V^{(l+1)^T}_{2i}\\
			U^{(l+1)}_{2i} V^{(l+1)^T}_{2i-1} &
			K_{2i}^{(l+1)} 
		\end{bmatrix}\label{eq:hodlr_gen}
	\end{equation}
	
	where $K_i^{(l)} \in \Rb^{N/2^l}$, $K_{2i-1}^{(l+1)}, K_{2i}^{(l+1)} \in \mathbb{R}^{N/2^{l+1} \times N/2^{l+1}}$, $U^{(l+1)}_{2i-1}, V^{(l+1)}_{2i}, U^{(l+1)}_{2i} V^{(l+1)}_{2i-1} \in \Rb^{N/2^{l+1} \times p}$, where $p \ll N$. Figure~\ref{fig:hodlr_mat} shows the HODLR representation of the dense matrix $A$ for different levels. Using this HODLR representation of the dense matrix A, the matrix-vector product can be performed in $\mathcal{O}(pN\log (N))$. The implementation of HODLR can be found in \cite{HODLRlib}.
	
	\section{\textbf{Rank growth of different interactions in \texorpdfstring{$2$}{2}D}}
	\label{sec3}
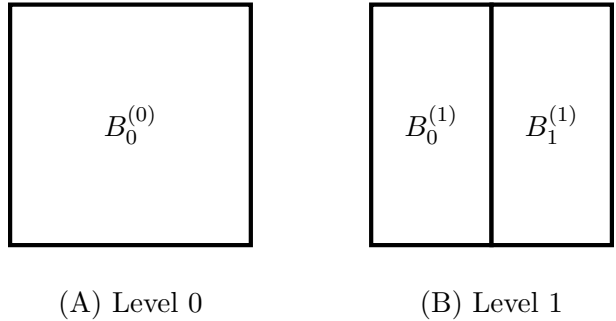
\begin{figure}[!htbp]
		\centering
		\begin{tikzpicture}[scale=0.8]
			\rect{white}{0,0}{4,4}
			\rect{white}{6,0}{8,4}
			\rect{white}{10,0}{8,4}
			\node at (2,2) {$B_0^{(0)}$};
			\node at (7,2) {$B_0^{(1)}$};
			\node at (9,2) {$B_1^{(1)}$};
			\node at (2,-1) {(A) Level $0$};
			\node at (8,-1) {(B) Level $1$};
		\end{tikzpicture}
		\caption{Subdivision of box $B_0^{(0)}$ into $B_0^{(1)}$ and  $B_1^{(1)}$}
		\label{bdomain01}
	\end{figure}

	In this section, we will examine the rank of the off-diagonal blocks of HODLR matrix arising out of $N$ body problem in two dimensions. Consider the box $B_0^{(0)} = [-1,1]^2$ containing $2n$ particles located at $\{\vec{r}_i\}_{i=1}^{2n}$, where $\vec{r}_i \in \Rb^2$. Let $K \in \Rb^{2n \times 2n}$ be the interaction matrix whose entries are given by
	\begin{equation}
		K(i,j) = 
		\begin{cases}
			0 & \text{ if }i =j\\
			\log\bkt{\abs{\vec{r}_i-\vec{r}_j}} & \text{ if }i \neq j
		\end{cases}
	\end{equation}
	Let's now represent $K$ as a $1$ level HOLDR matrix.
	As mentioned before, HODLR uses a K-D tree, and hence at each level in the hierarchical tree, the box is subdivided into two smaller boxes each containing $n$ particles as shown in Figure~\ref{bdomain01}. Hence, $B_0^{(0)}$ is subdivided into two boxes $B_0^{(1)} = [-1,0)\times [-1,1]$ and $B_1^{(1)} = [0,1]\times [-1,1]$ at level $1$. Let the boxes $B_0^{(1)},B_1^{(1)}$ each contain $n = m^2$ particles located at $m \times m$ Chebyshev grid, where $m \in \Zb^+$. Let $K_E\in \mathbb{R}^{n\times n}$ represent the interaction of particles in $B_1^{(1)}$ with the particles in $B_0^{(1)}$.
	% HIERARCHICAL SUBDIVISION
	
	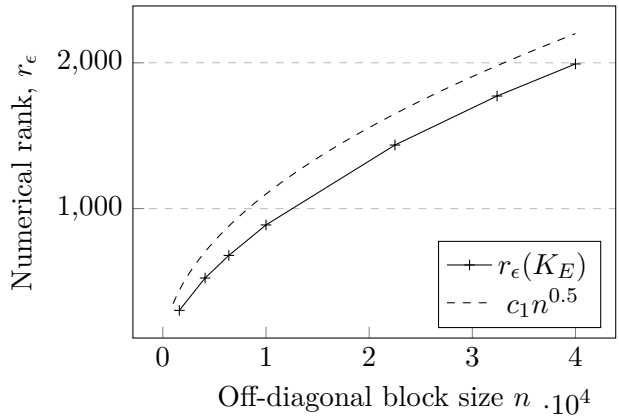
\begin{figure}[!htbp]
		\centering
		\begin{tikzpicture}
			\pgfplotsset{compat=newest, width=8cm, height=6cm,xtick pos=left, ytick pos=left}
			\begin{axis}[xlabel={Off-diagonal block size $n$},
				ylabel={Numerical rank, $r_\epsilon$},
				legend pos=south east,
				ymajorgrids=true,
				grid style=dashed,]
				% Edge Sharing
				\addplot[color=black,mark=+,] coordinates {
					(1600,302)
					(4096,524)
					(6400,679)
					(10000,888)
					(22500,1436)
					(32400,1772)
					(40000,1992)
				};
				\addplot [
				domain=1000:40000, 
				samples=500, 
				color=black,dashed,
				]
				{11*x^0.5};
				\legend{
					$r_{\epsilon}{(K_E)}$,
					$c_1n^{0.5}$
				}
			\end{axis}
		\end{tikzpicture}%
		\caption{Scaling of the rank of the off-diagonal block with the entries defined by $\log(r)$ to its size $n$}
		\label{fig:rank_edge}
	\end{figure}
	
	From Figure~\ref{fig:rank_edge}, we see that the numerical rank of the interaction matrix, $K_E$, increases as a function of $n$ and the scaling seems to be $\mathcal{O}\bkt{\sqrt{n}}$ (We prove the precise form of the scaling for uniformly distributed particles at the end of this section, Theorem~\ref{tm:farfield}). This immediately implies that if we were to use a HODLR matrix to represent the dense matrix $K$, the matrix-vector product will no longer scale linearly in $n$. This is because the rank of the largest off-diagonal block $K_E$ seems to scale as $\mathcal{O}\bkt{\sqrt{n}}$ (as shown in Figure~\ref{fig:rank_edge}). Hence, using HODLR one cannot obtain an almost linear scaling algorithm for $N$ body problems in $2$D.
	
	Consider the subdivision of the Box $B_0^{(0)}$ as in Figure \ref{fig:xv}, the Box marked as $X$ denotes a cluster with $n$ particles and let its respective index set be $I_X$. We consider three Boxes: $E$, $V$ and $F$, where Box $E$ shares an edge, $V$ shares a vertex with cluster $X$ respectively. Box $F$ does not share a boundary with Box $X$. The particles in $E$, $V$ and $F$, are indexed using the index sets $I_E$, $I_V$ and $I_F$, respectively. Let us examine the numerical rank of the interaction matrices $K\bkt{I_X,I_E}$, $K\bkt{I_X,I_V}$ and $K\bkt{I_X,I_F}$ (submatrices represented using MATLAB notation) for the kernel $\log(r)$.
	
	\begin{figure}[H]
		\centering
		\begin{tikzpicture}
			
			\draw[step=1cm,color=gray] (0,0) grid (4,4);
			
			\draw (0,0) rectangle (1,1) node[pos=0.5,text=black] {};
			\draw (1,0) rectangle (2,1) node[pos=0.5,text=black] {};
			\draw (0,1) rectangle (1,2) node[pos=0.5,text=black] {};
			\draw (1,1) rectangle (2,2) node[pos=0.5,text=black] {X};
			
			\draw (2,0) rectangle (3,1) node[pos=0.5,text=black] {};
			\draw (3,0) rectangle (4,1) node[pos=0.5,text=black] {};
			\draw (2,1) rectangle (3,2) node[pos=0.5,text=black] {E};
			\draw (3,1) rectangle (4,2) node[pos=0.5,text=black] {};
			
			\draw (0,2) rectangle (1,3) node[pos=0.5,text=black] {};
			\draw (1,2) rectangle (2,3) node[pos=0.5,text=black] {};
			\draw (0,3) rectangle (1,4) node[pos=0.5,text=black] {};
			\draw (1,3) rectangle (2,4) node[pos=0.5,text=black] {};
			
			\draw (2,2) rectangle (3,3) node[pos=0.5,text=black] {V};
			\draw (3,2) rectangle (4,3) node[pos=0.5,text=black] {};
			\draw (2,3) rectangle (3,4) node[pos=0.5,text=black] {};
			\draw (3,3) rectangle (4,4) node[pos=0.5,text=black] {F};

			\clip (0,0) rectangle (4 ,4);
		\end{tikzpicture}
		\caption{Clusters at 2 level on subdivision of box $B_0^{(0)}$}
		\label{fig:xv}
	\end{figure}
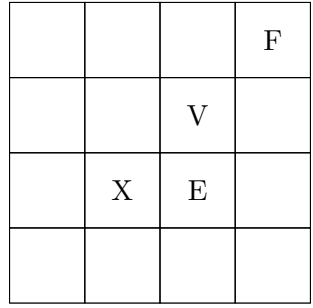
	\begin{figure}[H]
		\begin{subfigure}[b]{0.33\textwidth}
			\centering
			\includegraphics[width=\textwidth]{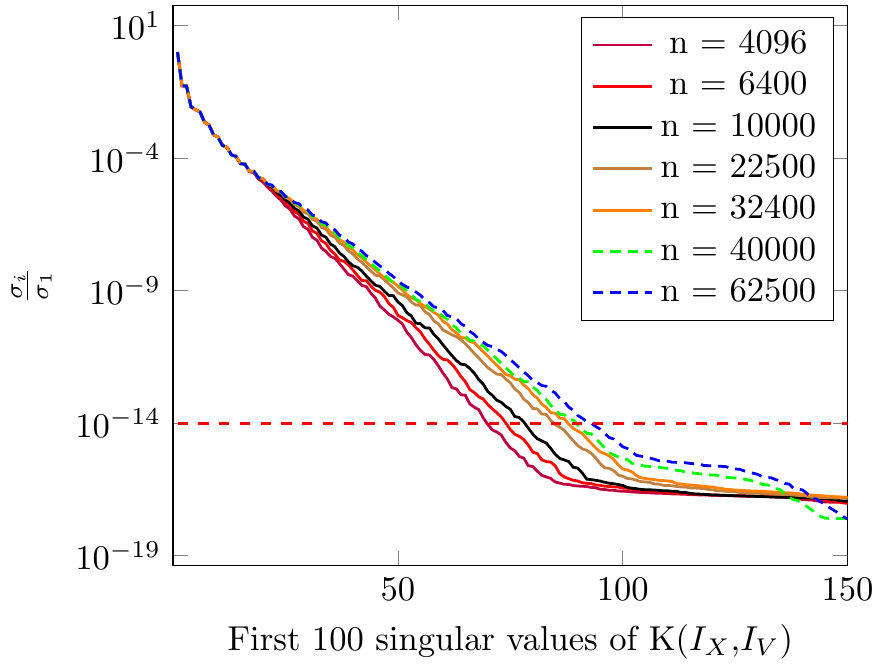}
			\caption{Vertex sharing clusters}
			\label{fig:l1}
		\end{subfigure}%
		\begin{subfigure}[b]{0.33\textwidth}
			\centering
			\includegraphics[width=\textwidth]{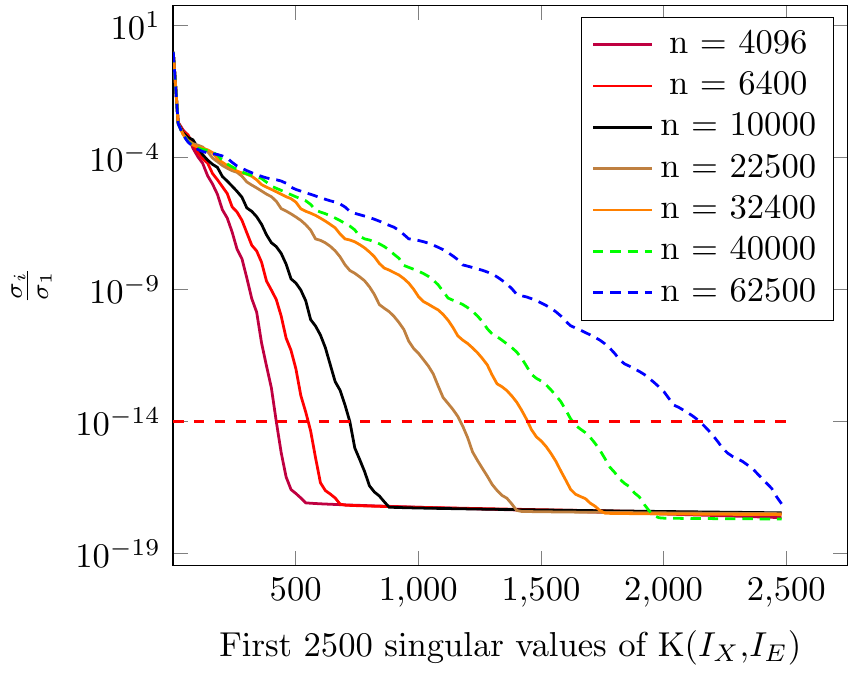}
			\caption{Edge sharing clusters}
			\label{fig:l2}
		\end{subfigure}%
		\begin{subfigure}[b]{0.33\textwidth}
			\centering
			\includegraphics[width=\textwidth]{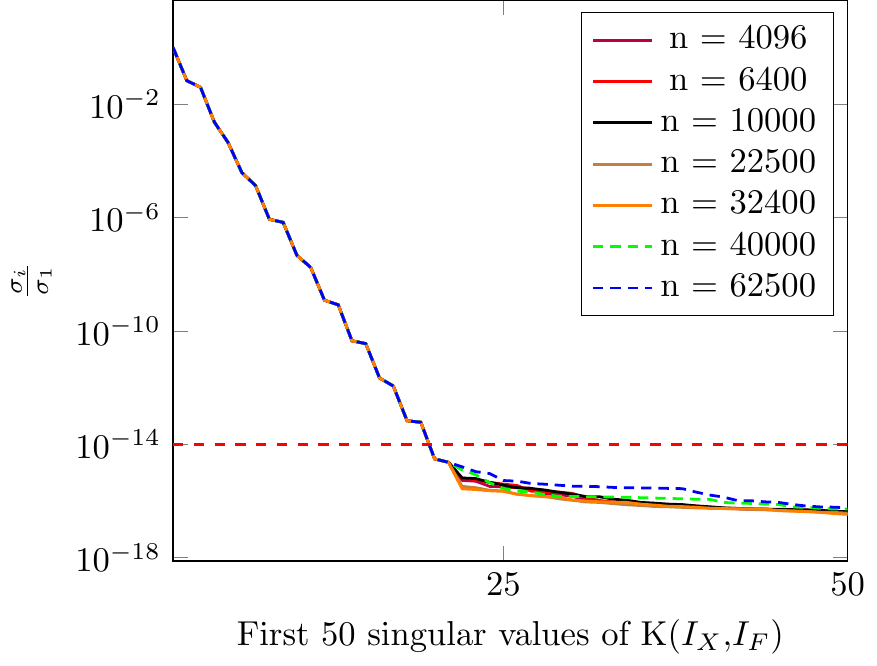}
			\caption{Far field clusters}
			\label{fig:ffs}
		\end{subfigure}%
		\caption{Decay of singular values of different submatrices of size ($n$) for the $\log (r)$ kernel}
		\label{fig:three_sing}
	\end{figure}

	From Figure~\ref{fig:three_sing}, we observe that the singular values of the matrix $K(I_X,I_F)$ decay rapidly, followed by the singular values of the matrix $K(I_X,I_V)$. The decay of singular values of the matrix $K(I_X,I_E)$ is the slowest among the three. We now prove theorems, which validate the observation in Figure~\ref{fig:three_sing}.
	
	\begin{lemma}
		\label{lemma}
		\noindent\textbf{(Multipole Expansion)} Consider $n$ charges of strength $\{q_i\}_{i=1}^n$ located at points $\{z_i\}_{i=1}^n$, where $z_i \in \mathbb{C}$ with $\abs{z_i} < r$. Then for any $z \in \mathbb{C}$ with $\abs{z} > r$, the complex potential induced by the $n$ charges is given by $\phi(z) = \dsum_{i=1}^n q_i \log\bkt{z-z_i}$. Note that the real potential, $\dsum_{i=1}^n q_i \log\bkt{\abs{z-z_i}}$, is nothing but the real part of $\phi(z)$, i.e., $\text{Re}\bkt{\phi(z)}$.\\
		We then have $$\phi(z) = q \log(z) + \dsum_{k=1}^{\infty} \dfrac{a_k}{z^k}$$
		where $q = \dsum_{i=1}^n q_i$ and $a_k = -\dsum_{i=1}^n \dfrac{q_iz_i^k}{k}$. Furthermore, for any $p \geq 1$,
		\begin{equation}
			\abs{\phi(z) - q\log(z)-\dsum_{k=1}^p \dfrac{a_k}{z^k}} \leq \bkt{\dfrac{Q}{p+1}} \bkt{\dfrac1{c-1}} \bkt{\dfrac1c}^p
		\end{equation}
		where $c=\abs{\dfrac{z}r}$, $Q = \dsum_{i=1}^m \abs{q_i}$
	\end{lemma}
	The proof for this lemma can be found in \cite{FMMcourse,Greengard_1987}.
	
	\begin{theorem}[Rank of different interactions in $2$D]
		\label{tm:farfield}
		Let $N$ charges, $\{q_j\}_{j=1}^N$, be uniformly located at $\{z_j\}_{j=1}^N$ inside a box $B_1$ of side length $a$. Let $Q = \dsum_{i=1}^n \abs{q_i}$. The complex potential due to these $N$ charges, at $M$ locations, $\{w_i\}_{i=1}^M$, inside another box $B_2$ is given by $\phi_i = \dsum_{j=1}^N \log\bkt{w_i-z_j}q_j$. In matrix-vector parlance, we have
		$$\vec{\phi} = A \vec{q}$$
		where $\vec{q} \in \Rb^{N \times 1}$, $\vec{\phi} \in \Cb^{M \times 1}$ and $A \in \Cb^{M \times N}$ with $A_{ij} = \log\bkt{w_i-z_j}$. Then we claim that given $\epsilon > 0$, there exists a matrix $\tilde{A} \in \Cb^{M \times N}$ with rank at most $p$, where $p \in \mathcal{O} \bkt{R(N)\log\bkt{R(N)Q/\epsilon}}$ such that $\abs{\phi_i-\bkt{\tilde{A}q}_i} < \epsilon$, where
		\begin{enumerate}[(i)]
			\item $R(N) = 1$ if the boxes $B_1$ and $B_2$ are one box away.
			\item $R(N) = \log(N)$ if the boxes $B_1$ and $B_2$ are vertex sharing neighbors.
			\item $R(N) = \sqrt{N}$ if the boxes $B_1$ and $B_2$ are edge sharing neighbors.
		\end{enumerate}
		\begin{proof}
			\begin{enumerate}[(i)]
				% Far field boxed
				\item Boxes $B_1$ and $B_2$ are one box away as shown in Figure~\ref{one_box_away}.
				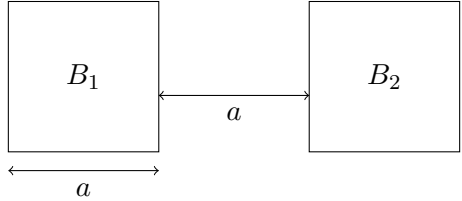
\begin{figure}[!htbp]
					\begin{center}
						\begin{tikzpicture}[scale=1]
							\draw (-1,-1) rectangle (1,1);
							\draw (3,-1) rectangle (5,1);
							\node at (0,0) {$B_1$};
							\node at (4,0) {$B_2$};
							\draw [<->] (-1,-1.25) -- (1,-1.25);
							\draw [<->] (1,-0.25) -- (3,-0.25);
							\node at (2,-0.5) {$a$};
							\node at (0,-1.5) {$a$};
						\end{tikzpicture}
					\end{center}
					\caption{Rank of Far field boxes}
					\label{one_box_away}
				\end{figure}
				The proof follows immediately from Lemma~\ref{lemma}. We construct a circle, say $C$, by enclosing the box $B_1$ of radius $\dfrac{a}{\sqrt2}$. The distance between the center of the circle $C$ and the box $B_2$ is $\dfrac{3a}2$. Hence in the Lemma~\ref{lemma}, we have $c = \dfrac{3a/2}{a/\sqrt2} = \dfrac3{\sqrt2}$. Now if we set $$\tilde{A}_{ij} = \log\bkt{w_i} - \dsum_{k=1}^p \dfrac{z_j^k}{kw_i^k}$$ we have
				$$\bkt{\tilde{A}q}_i = q\log{w_i} - \dsum_{j=1}^n \dsum_{k=1}^p \dfrac{q_jz_j^k}{kw_i^k}$$
				From here we see that
				$$\tilde{A}_{ij} = \begin{bmatrix}\log\bkt{w_i} & \dfrac{-1}{w_i} & \dfrac{-1}{2w_i^2} & \cdots &  \dfrac{-1}{pw_i^p}\end{bmatrix}
				\begin{bmatrix} 1 \\ z_j \\ z_j^2 \\ \vdots \\ z_j^p\end{bmatrix}$$
				Hence, we note that the rank of $\tilde{A}$ is $p+1$. Hence, from Lemma~\ref{lemma}, we have
				$$\abs{\bkt{Aq}_i-\bkt{\tilde{A}q}_i} \leq \dfrac{Q}{p+1} \bkt{\dfrac{\sqrt2}{3-\sqrt2}} \bkt{\dfrac{\sqrt2}{3}}^p$$
				Now choosing $p = \dfrac{\log\bkt{\sqrt2 Q/ \bkt{\bkt{3-\sqrt2}\epsilon}}}{\log{\bkt{3/\sqrt2}}}$ guarantees
				$$\abs{\bkt{Aq}_i-\bkt{\tilde{A}q}_i} \leq \epsilon $$
				So in this case, when the boxes are one box away, we have that the rank of the matrix $\tilde{A}$ to be $p+1$, where $\boxed{p+1 \in \mathcal{O} \bkt{\log\bkt{Q/\epsilon}}}$.
				% Vertex sharing boxes
				\item
				In this part, we study the growth in the rank of the interaction matrix corresponding to vertex sharing boxes as shown in Figure~\ref{vertex_sharing_box}.
				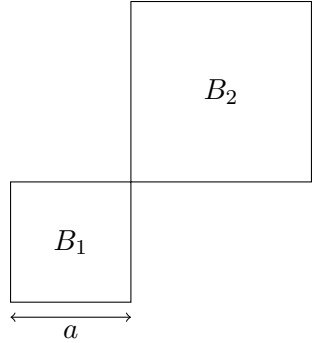
\begin{figure}[!htbp]
					\begin{center}
						\begin{tikzpicture}[scale=0.8]
							\draw (-1,-1) rectangle (1,1);
							\draw (1,1) rectangle (4,4);
							\node at (0,0) {$B_1$};
							\node at (2.5,2.5) {$B_2$};
							\draw [<->] (-1,-1.25) -- (1,-1.25);
							\node at (0,-1.5) {$a$};
						\end{tikzpicture}
					\end{center}
					\caption{Rank of vertex sharing boxes}
					\label{vertex_sharing_box}
				\end{figure}
				The proof for this part relies on hierarchically subdividing the box $B_1$ as shown in Figure~\ref{hierarchical_division_box}.
				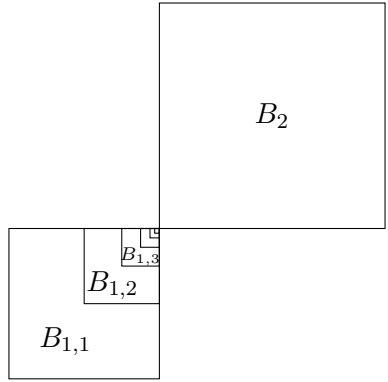
\begin{figure}
					\centering
					\begin{tikzpicture}
						\draw (-1,-1) rectangle (1,1);
						\draw (0,0) rectangle (1,1);
						\draw (0.5,0.5) rectangle (1,1);
						\draw (0.75,0.75) rectangle (1,1);
						\draw (0.875,0.875) rectangle (1,1);
						\draw (0.9375,0.9375) rectangle (1,1);
						\draw (1,1) rectangle (4,4);
						\node at (-0.25,-0.5) {$B_{1,1}$};
						\node at (0.375,0.25) {$B_{1,2}$};
						\node at (0.75,0.625) {\tiny $B_{1,3}$};
						\node at (2.5,2.5) {$B_2$};
					\end{tikzpicture}
					\caption{Hierarchical subdivision of box $B_1$}
					\label{hierarchical_division_box}
				\end{figure}
				$B_{1,j}$ denotes the L-shaped domain at the $j^{th}$ level and $B_{1,j}$ partitions the box $B_1$, i.e., $B_{1,j} \dcap B_{1,k} = \emptyset$ and $B_1 = \displaystyle \bigcup_{j=1}^{\infty} B_{1,j}$. Note that since we have $N$ particles uniformly distributed in the domain, there exists $\kappa=\log_4\bkt{N}+\text{constant}$, beyond which there will be no particles in the L-shaped box, i.e., there will be no particles in the boxes $B_{1,j}$ for $j >\kappa$. Let $A_j \in \Cb^{M \times N}$ be the matrix whose columns corresponding to charges lying inside $B_{1,j}$ are the respective columns of the matrix $A$ and the columns corresponding to the charges lying outside $B_{1,j}$ are zeroes.
			
				We will make use of the multipole expansion lemma (Lemma~\ref{lemma}) to prove the above claim. To do that, we will construct a circle enclosing $B_{1,1}$ as shown in Figure~\ref{fig_3}.
				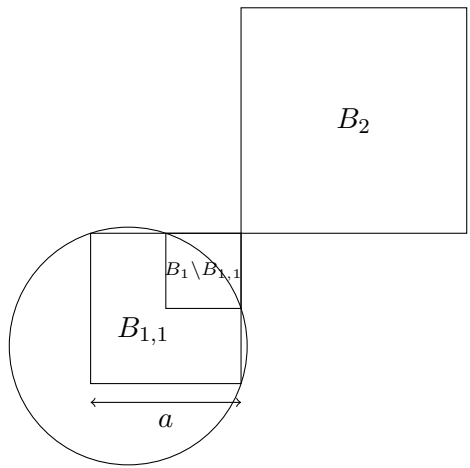
\begin{figure}[!htbp]
					\begin{center}
						\begin{tikzpicture}
							\draw (-1,-1) rectangle (1,1);
							\draw (0,0) rectangle (1,1);
							\draw (1,1) rectangle (4,4);
							\node at (-0.3,-0.3) {$B_{1,1}$};
							\node at (0.5,0.5) {\tiny $B_1\backslash B_{1,1}$};
							\node at (2.5,2.5) {$B_2$};
							\draw [<->] (-1,-1.25) -- (1,-1.25);
							\node at (0,-1.5) {$a$};
							\draw (-0.5,-0.5) circle (1.58113);
						\end{tikzpicture}
					\end{center}
					\caption{Rank of sub-boxes}
					\label{fig_3}
				\end{figure}
				The radius of the circle can be shown to be $r= \dfrac{a\sqrt{10}}4$ and the center to be at $(a/4,a/4)$ (assuming the box $B_1 = [0,a]^2$). The locations at which the potential is measured is bounded below by $\abs{z} > \dfrac{3a\sqrt2}4$ (essentially the distance of the top right of box $B_1$ from the center of the circle). Hence, for the Multipole expansion lemma, we have $c = \abs{\dfrac{z}r} > \dfrac3{\sqrt5}$. Hence, we obtain
				$$\abs{\bkt{A_{1}q-\tilde{A}_{1}q}_i} \leq \bkt{\dfrac{Q}{p+1}} \bkt{\dfrac{\sqrt5}{3-\sqrt5}} \bkt{\dfrac{\sqrt5}3}^{p}$$
				where $Q = \dsum_{i=1}^m \abs{q_i}$ and $\tilde{A}_1 \in \Cb^{M \times N}$ is a matrix of rank $p+1$ obtained from the multipole expansion as in the previous proof (Note that $p$ is independent of the box $B_{1,j}$). Note that, as with the matrix $A_1$, the columns corresponding to the charges lying outside $B_{1,1}$ are zeroes. Now choosing $p = \ceil{\dfrac{\log_4 \bkt{\dfrac{4Q}{\epsilon_1}}}{\log_4(3/\sqrt5)}}$ guarantees
				$$\abs{\bkt{A_{1}q-\tilde{A}_{1}q}_i} <  \epsilon_1$$
				Now repeat the same for the box $B_1 \backslash B_{1,1}$, i.e., split the charges in $B_1 \backslash B_{1,1}$ as those in the domain $B_{1,2}$ and those outside the domain $B_{1,2}$. Let $A_{2} \in \Cb^{M \times N}$ be the matrix that corresponds to charges in the domain $B_{1,2}$. By a similar argument as above, we have that the matrix $A_{2}$ can be approximated by a matrix $\tilde{A}_{2} \in \Cb^{M \times N}$ of rank $1+\ceil{\dfrac{\log_4\bkt{\dfrac{4Q}{\epsilon_1}}}{\log_4(3/\sqrt5)}}$ such that $\abs{\bkt{A_{2}q-\tilde{A}_{2}q}_i} < \epsilon_1$.\\
				Repeating this till $\kappa$ levels, we have the rank of the matrix $\tilde{A} = \tilde{A}_{1} + \tilde{A}_{2} + \cdots + \tilde{A}_{\kappa}$, where $\tilde{A}_k,\tilde{A} \in \Cb^{M \times N}$, to be bounded above by $\bkt{1+\ceil{\dfrac{\log_4\bkt{\dfrac{4Q}{\epsilon_1}}}{\log_4(3/\sqrt5)}}}\kappa$. The matrix $A$ is now approximated by
				$\tilde{A} = \tilde{A}_{1} + \tilde{A}_{2} + \cdots + \tilde{A}_{\kappa}$
				whose rank is bounded above $\bkt{1+\ceil{\dfrac{\log_4\bkt{\dfrac{4Q}{\epsilon_1}}}{\log_4(3/\sqrt5)}}}\kappa$ and the error is bounded above by $\kappa\epsilon_1$. This is because we have
				\begin{align*}
					\abs{\bkt{Aq-\tilde{A}q}_i} & = \abs{\bkt{\bkt{A_1+A_2+\cdots+A_{\kappa}}q-\bkt{\tilde{A}_1+\tilde{A}_2+\cdots+\tilde{A}_{\kappa}}q}_i} \\
					& \leq \abs{\bkt{A_1q-\tilde{A}_1q}_i} + \abs{\bkt{A_2q-\tilde{A}_2q}_i} + \cdots + \abs{\bkt{A_{\kappa}q-\tilde{A}_{\kappa}q}_i}\\
					& < \kappa \epsilon_1
				\end{align*}
				$$$$
				Hence, to obtain $\tilde{A}$ such that $\abs{\phi_i-\bkt{\tilde{A}q}_i} < \epsilon$, we need to pick $\epsilon_1 = \dfrac{\epsilon}{\kappa}$.\\
				Hence, given $\epsilon > 0$, there exists a matrix $\tilde{A}$ with rank at most
				$$\bkt{1+\ceil{\dfrac{\log_4\bkt{\dfrac{4Q\kappa}{\epsilon}}}{\log_4(3/\sqrt5)}}}\kappa$$
				such that $\abs{\phi_i-\bkt{\tilde{A}q}_i} < \epsilon$.\\
				Hence, we have the rank in this case to be $\mathcal{O}\bkt{\log\bkt{N}\log\bkt{\dfrac{Q\log\bkt{N}}{\epsilon}}}$, since $\kappa = \log_4\bkt{N}+\text{constant}$.
				
				\item
				In the final part, we study the growth in rank of edge sharing boxes as shown in Figure~\ref{edge_sharing_box}.
				\begin{figure}[!htbp]
					\begin{center}
						\begin{tikzpicture}[scale=1]
							\draw (-1,-1) rectangle (1,1);
							\draw (1,-1) rectangle (3,1);
							\node at (0,0) {$B_1$};
							\node at (2,0) {$B_2$};
							\draw [<->] (-1,-1.25) -- (1,-1.25);
							\node at (0,-1.5) {$a$};
						\end{tikzpicture}
					\end{center}
					\caption{Rank of edge sharing boxes}
					\label{edge_sharing_box}
				\end{figure}
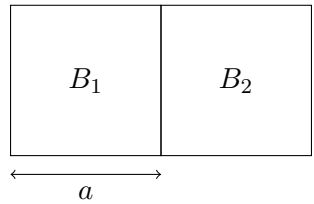
			\end{enumerate}
			The proof again relies on hierarchically subdividing the box $B_1$ like the earlier proof for vertex sharing interactions but the hierarchical subdivision is done in a different way as shown in Figure~\ref{hierarchical_division_box_edge}.
			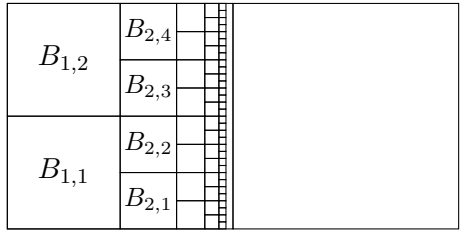
\begin{figure}[!htbp]
				\centering
				\begin{tikzpicture}[scale=1.5]
					\draw (-1,-1) rectangle (1,1);
					\draw (1,-1) rectangle (3,1);
					\draw (-1,-1) rectangle (0,0);
					\draw (-1,0) rectangle (0,1);
					
					\foreach \i in {0,1,2,3}
					\draw (0,-1+0.5*\i) rectangle (0.5,-0.5+0.5*\i);
					
					\foreach \i in {0,1,2,...,7}
					\draw (0.5,-1+0.25*\i) rectangle (0.75,-0.75+0.25*\i);
					
					\foreach \i in {0,1,2,...,15}
					\draw (0.75,-1+0.125*\i) rectangle (0.875,-0.875+0.125*\i);
					
					\foreach \i in {0,1,2,...,31}
					\draw (0.875,-1+0.0625*\i) rectangle (0.9375,-0.9375+0.0625*\i);
					
					\node at (-0.5,-0.5) {$B_{1,1}$};
					\node at (-0.5,0.5) {$B_{1,2}$};
					
					\node at (0.25,-0.75) {\small $B_{2,1}$};
					\node at (0.25,-0.25) {\small $B_{2,2}$};
					\node at (0.25,0.25) {\small $B_{2,3}$};
					\node at (0.25,0.75) {\small $B_{2,4}$};
					
				\end{tikzpicture}
				\caption{Hierarchical subdivision of box $B_1$}
				\label{hierarchical_division_box_edge}
			\end{figure}
			The box $B_{k,i}$ denotes the $i^{th}$ box from the bottom at the $k^{th}$ level in the hierarchical subdivision, where $1 \leq k < \infty$ and $1 \leq i \leq 2^k$. Note that $B_{k,i}$ partitions the box $B_1$ (i.e., $B_{k,r} \dcap B_{j,q}$ is non-empty iff $k=j$ and $r=q$; $B_1 = \dbcup_{k=1}^{\infty} \dbcup_{i=1}^{2^k} B_{k,i}$). Note that the boxes beyond a certain level $\kappa = \log_4\bkt{N}+\text{constant}$ will no longer contain any charges (since we have assumed uniform distribution of charges in the domain). Let $A_{k,i} \in \Cb^{M \times N}$ be the matrix whose columns corresponding to charges lying inside $B_{k,i}$ are the respective columns of the matrix $A$ and the columns corresponding to the charges lying outside $B_{k,i}$ are zeroes.
			
			We will again rely on the multipole expansions (Lemma~\ref{lemma}) to prove this claim. We construct a circle, say $C_{k,j}$, enclosing the box $B_{k,j}$. The radius of the circle is $r=\dfrac{a}{2^k\sqrt2}$. The distance of the box $B_2$ from the center of the circle $C_{k,i}$ (which is same as the center of the box $B_{k,i}$) is $\dfrac{3a}{2^{k+1}}$. Hence, we have $c = \abs{\dfrac{z}{r}} > \dfrac{3a/2^{k+1}}{a/2^{k+1/2}} = \dfrac{3}{\sqrt2}$. Hence, we obtain
			$$\abs{\bkt{A_{k,j}q-\tilde{A}_{k,j}q}_i} \leq \dfrac{Q}{p+1} \bkt{\dfrac{\sqrt2}{3-\sqrt2}} \bkt{\dfrac{\sqrt2}3}^{p}$$
			where $Q = \dsum_{i=1}^m \abs{q_i}$ and $\tilde{A}_{k,j} \in \Cb^{M \times N}$ is a matrix of rank $p+1$ obtained from the multipole expansion as in the proof of the first part. Note that, as with the matrix $A_{k,j}$, the columns corresponding to the charges lying outside $B_{k,j}$ are zeroes. Now choosing $p = \ceil{\dfrac{\log \bkt{Q/\epsilon_1}}{\log \bkt{3/\sqrt2}}}$ (Note that $p$ doesn't depend on the box level or box number) guarantees
			$$\abs{\bkt{A_{k,j}q-\tilde{A}_{k,j}q}_i} < \epsilon_1$$
			Repeating the same for all boxes till level $\kappa$, we have that the rank of the matrix $\tilde{A} = \dsum_{k=1}^{\kappa} \dsum_{j=1}^{2^k} \tilde{A}_{k,j}$, where $\tilde{A}_{k,j}, \tilde{A} \in \Cb^{M \times N}$, to be bounded above by
			$$\dsum_{k=1}^{\kappa} 2^k \bkt{1+\ceil{\dfrac{\log \bkt{Q/\epsilon_1}}{\log \bkt{3/\sqrt2}}}} = \bkt{2^{\kappa+1}-2} \bkt{1+\ceil{\dfrac{\log \bkt{Q/\epsilon_1}}{\log \bkt{3/\sqrt2}}}} \in \mathcal{O}\bkt{\sqrt{N}\log\bkt{Q/\epsilon_1}}$$ and the error is bounded above by
			$$\dsum_{k=1}^{\kappa} 2^k \epsilon_1 = \bkt{2^{\kappa+1}-2}\epsilon_1 <2\sqrt{N}\epsilon_1$$
			Hence, to obtain $\tilde{A}$ such that
			$$\abs{\phi_i-\bkt{\tilde{A}q}_i} < \epsilon$$
			we need to pick $\epsilon_1 = \dfrac{\epsilon}{2\sqrt{N}}$.\\
			Hence, given $\epsilon > 0$, there exists a matrix $\tilde{A}$ with rank $p \in \mathcal{O}\bkt{\sqrt{N} \log\bkt{\sqrt{N}Q/\epsilon}}$ such that
			$$\abs{\phi_i - \bkt{\tilde{A}q}_i} < \epsilon$$
		\end{proof}
	\end{theorem}
	
	\begin{figure}[H]
		\begin{subfigure}[b]{0.45\textwidth}
			\centering
			\includegraphics[width=\textwidth]{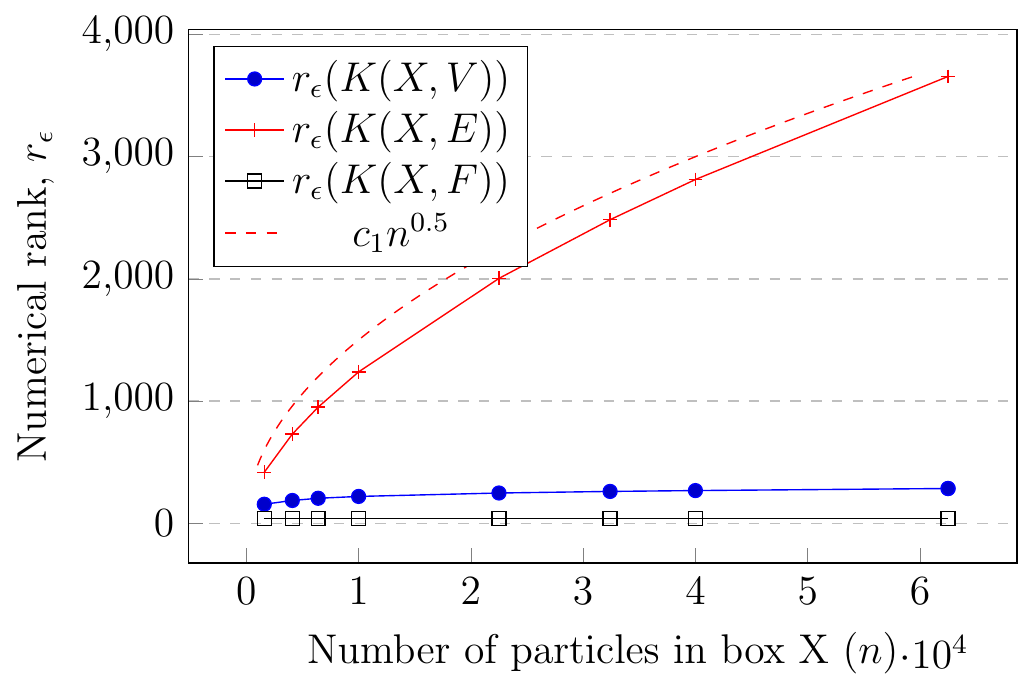}
			\caption{$\frac{1}{r}$}
			\label{fig:lrk1}
		\end{subfigure}%
		\hfill
		\begin{subfigure}[b]{0.45\textwidth}
			\centering
			\includegraphics[width=\textwidth]{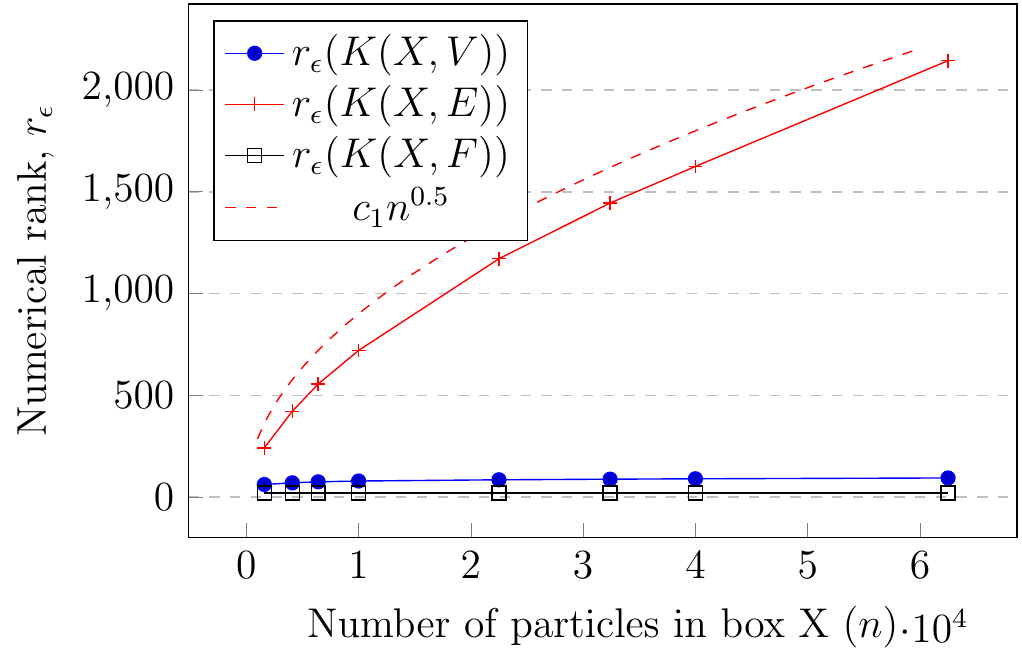}
			\caption{$\log(r)$}
			\label{fig:lrk2}
		\end{subfigure}%
		\hfill
		\begin{subfigure}[b]{0.45\textwidth}
			\centering
			\includegraphics[width=\textwidth]{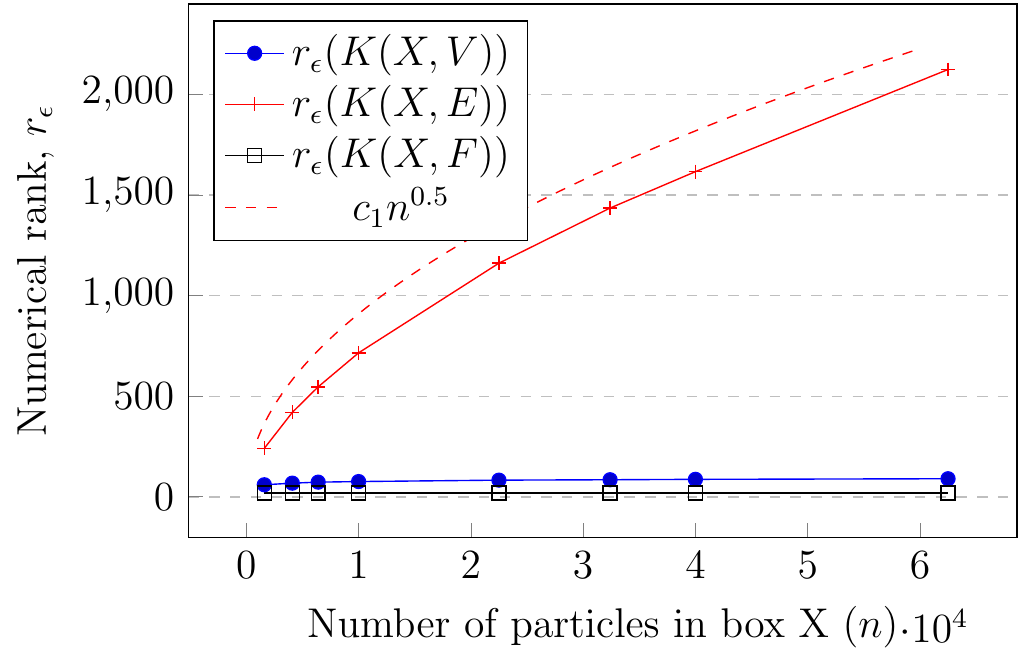}
			\caption{Hankel Function of first kind order zero}
			\label{fig:lrk3}
		\end{subfigure}%
		\hfill
		\begin{subfigure}[b]{0.45\textwidth}
			\centering
			\includegraphics[width=\textwidth]{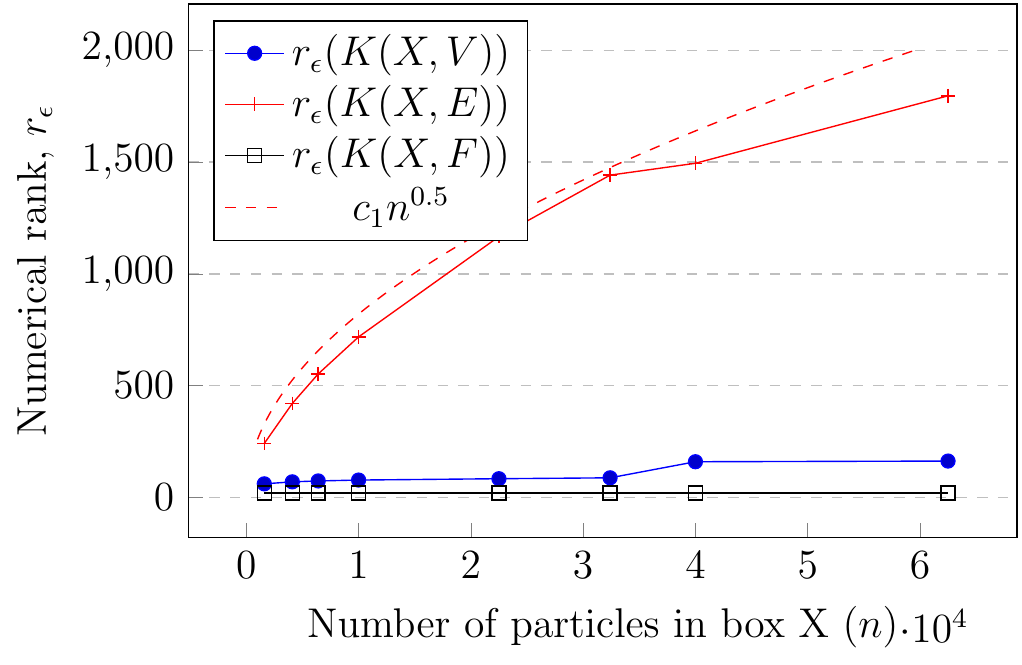}
			\caption{Bessel Function of second kind order zero}
			\label{fig:lrk4}
		\end{subfigure}%
		\hfill
		\begin{subfigure}[b]{0.45\textwidth}
			\centering
			\includegraphics[width=\textwidth]{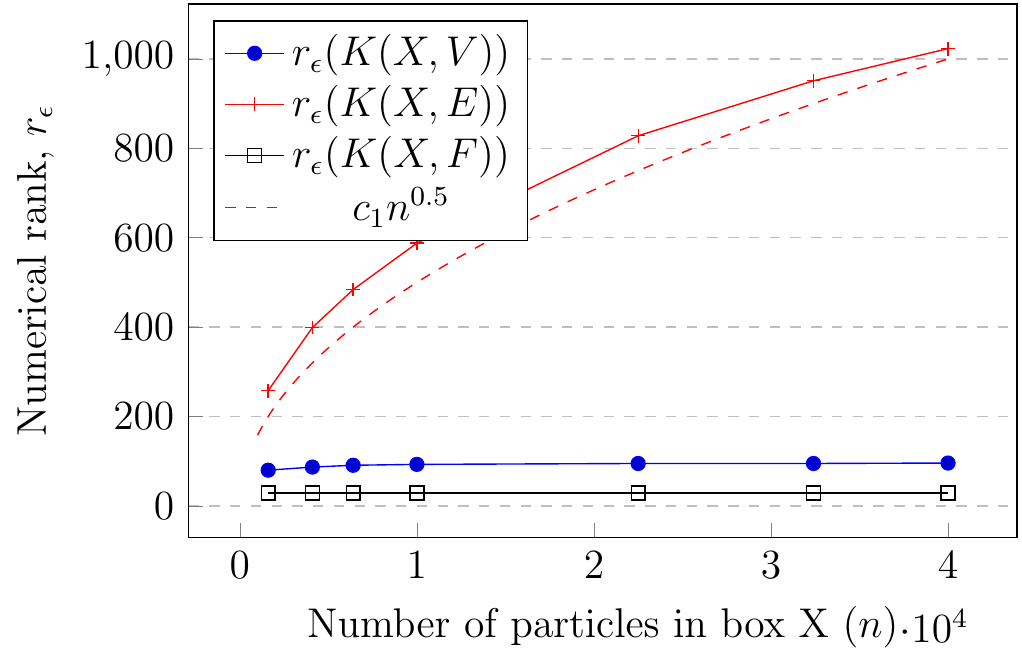}
			\caption{Thin Plate Spline; $r^2\log(r)$}
			\label{fig:lrk8}
		\end{subfigure}%
		\caption{Numerical Rank with $\epsilon=10^{-14}$   of different interactions vs size of the off-diagonal block (n)}
		\label{fig:three graphs}
	\end{figure}
	
	Figure~\ref{fig:three graphs} shows the numerical rank of the matrix corresponding to the three different interactions for a wide range of kernels arising from radial basis function interpolation and integral equations. It is evident from Figure~\ref{fig:three graphs} that the numerical rank remains constant for the far-field interaction, i.e., the numerical rank of the matrix $K(I_X,I_F)$. The numerical rank remains "almost" constant for the vertex sharing interaction, i.e., the numerical rank of the matrix $K(I_X,I_V)$, whereas the numerical rank increases with $n$ (roughly as $\sqrt{N}$) for the edge sharing interaction, i.e., the numerical rank of the matrix $K(I_X,I_E)$. The theoretical bounds from the Theorem~\ref{tm:farfield} is in agreement with the numerical experiments performed and reported in Figure~\ref{fig:three graphs}. We leverage these results to develop our new hierarchical low-rank structure, HODLR2D.
	\section{HODLR2D}
\label{sec4}
This section provides a detailed description of HODLR2D, the new hierarchical low-rank matrix structure. In the previous section (Section~\ref{sec3}), we provided numerical illustrations and a proof for the fact that the numerical rank of the matrix corresponding to edge sharing interactions grew as $\mathcal{O}\bkt{\sqrt{N}\log\bkt{\dfrac{\sqrt{N} Q}{\epsilon}}}$, where $N$ is the size of the off-diagonal submatrix. This immediately implies that the HODLR~\cite{ambikasaran2013large} based matrix-vector products will no longer scale linearly.

Further, in the previous section (Section~\ref{sec3}), we provided numerical illustrations and a proof for the fact that the numerical rank of the matrix corresponding to vertex sharing interactions grew as $\mathcal{O}\bkt{\log(N)\log\bkt{\dfrac{\log\bkt{\log(N) Q}}{\epsilon}}}$, where $N$ is the corresponding matrix size. We leverage this slow growth in the numerical rank of the matrix corresponding to vertex sharing interactions to construct the new HODLR2D structure.

\begin{figure}[H]
	\centering
	\begin{subfigure}{0.325\textwidth}
		\centering
		\begin{tikzpicture}[scale=1]
			%\coordinate[label=left:0]   (a) at (0,0);
			%	\coordinate[label=right:1]   (b) at (4,0);
			%	\coordinate[label=left:1]   (c) at (0,4);
			\draw [ultra thick](0,0) rectangle (4 ,4) node[pos=0.5,text=black] {0};
			\clip (0,0) rectangle (4 ,4);
			\pgfmathsetseed{24122015}
			\foreach \p in {1,...,1000}
			{
				\pgfmathsetmacro{\x}{4*rand}
				\pgfmathsetmacro{\y}{4*rand}
				\fill[red]    (\x,\y) circle (0.01);
			}
		\end{tikzpicture}
		\caption{Level $0$}
		\label{dlv0}
	\end{subfigure}
	\begin{subfigure}{0.325\textwidth}
		\centering
		\begin{tikzpicture}[scale=1]
			\coordinate[label=0]   (d) at (1,1);
			\coordinate[label=1]   (d) at (3,1);
			\coordinate[label=3]   (d) at (1,3);
			\coordinate[label=2]   (d) at (3,3);
			
			\draw[step=2cm,ultra thick] (0,0) grid (4,4);
			
			\clip (0,0) rectangle (4 ,4);
			\pgfmathsetseed{24122015}
			\foreach \p in {1,...,1000}
			{
				\pgfmathsetmacro{\x}{4*rand}
				\pgfmathsetmacro{\y}{4*rand}
				\fill[red]    (\x,\y) circle (0.01);
			}
		\end{tikzpicture}
		\caption{Level $1$}
		\label{dlv1}
	\end{subfigure}
	\begin{subfigure}{0.325\textwidth}
		\centering
		\begin{tikzpicture}[scale=1]
			\draw[step=1cm,ultra thick] (0,0) grid (4,4);
			
			\draw (0,0) rectangle (1,1) node[pos=0.5,text=black] {0};
			\draw (1,0) rectangle (2,1) node[pos=0.5,text=black] {1};
			\draw (0,1) rectangle (1,2) node[pos=0.5,text=black] {3};
			\draw (1,1) rectangle (2,2) node[pos=0.5,text=black] {2};
			
			\draw (2,0) rectangle (3,1) node[pos=0.5,text=black] {4};
			\draw (3,0) rectangle (4,1) node[pos=0.5,text=black] {5};
			\draw (2,1) rectangle (3,2) node[pos=0.5,text=black] {7};
			\draw (3,1) rectangle (4,2) node[pos=0.5,text=black] {6};
			
			\draw (0,2) rectangle (1,3) node[pos=0.5,text=black] {12};
			\draw (1,2) rectangle (2,3) node[pos=0.5,text=black] {13};
			\draw (0,3) rectangle (1,4) node[pos=0.5,text=black] {15};
			\draw (1,3) rectangle (2,4) node[pos=0.5,text=black] {14};
			
			\draw (2,2) rectangle (3,3) node[pos=0.5,text=black] {8};
			\draw (3,2) rectangle (4,3) node[pos=0.5,text=black] {9};
			\draw (2,3) rectangle (3,4) node[pos=0.5,text=black] {11};
			\draw (3,3) rectangle (4,4) node[pos=0.5,text=black] {10};

			\clip (0,0) rectangle (4 ,4);
			\pgfmathsetseed{24122015}
			\foreach \p in {1,...,1000}
			{
				\pgfmathsetmacro{\x}{4*rand}
				\pgfmathsetmacro{\y}{4*rand}
				\fill[red]    (\x,\y) circle (0.01);
			}
		\end{tikzpicture}
		\caption{Level $2$}
		\label{dlv2}
	\end{subfigure}
	\caption{Subdivision of the box $B\in \Rb^2$ at different levels and the numbering convention used}\label{sdomain01}
\end{figure}
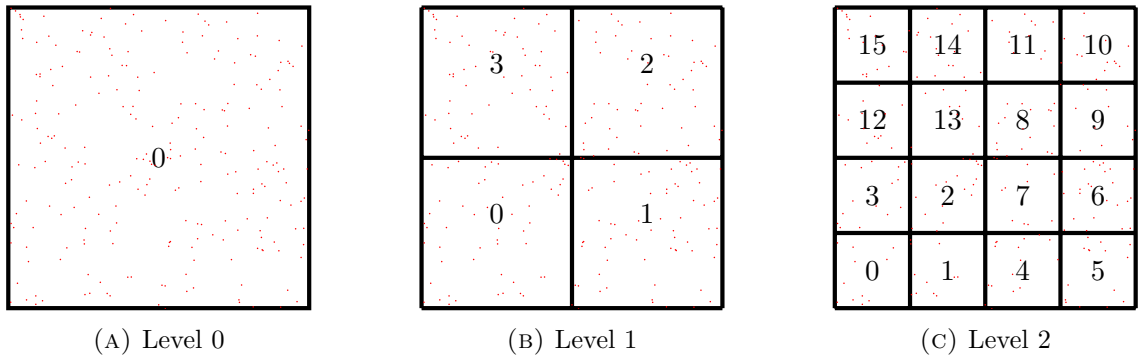

We subdivide the domain using a quadtree. We begin by subdividing the square box B (the entire domain) into $4$ smaller square boxes and number them as in Figure~\ref{sdomain01} and continue this till $\kappa$ levels. 
% Division of boxes

\begin{figure}[H]
	\centering
	\includegraphics[width=0.5\linewidth]{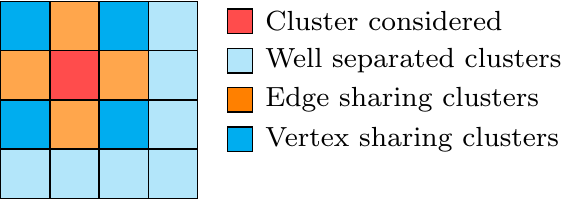}
	\caption{Types of clusters at level 2 for a cluster considered}
	\label{fig:clusRel}
\end{figure}

For each smaller box, the following definitions are intended to help us describe algorithms for HODLR2D matrices. Let $\mathcal{C}$ denote the cluster corresponding to the node $i$ at level $l$, $\mathcal{N}^{(l)}_i$ as defined in Section~\ref{sec2}.
\begin{definition}
	\label{def:ec}
	Two clusters in a level are said to be \textbf{edge sharing clusters}, if their respective boxes share an edge. Then, $\mathcal{E}_{\mathcal{C}}$ for a cluster $\mathcal{C}$ be a set containing the clusters that share an edge with $\mathcal{C}$. Note that $\abs{\mathcal{E}_{\mathcal{C}}} \leq 4$.
\end{definition}
\begin{definition}
	\label{def:vc}
	Two clusters in a level are said to be \textbf{vertex sharing clusters} if their respective boxes share a vertex. Then, $\mathcal{V}_{\mathcal{C}}$ for a cluster $\mathcal{C}$ be a set containing clusters that share a vertex with $\mathcal{C}$. Note that $\abs{\mathcal{V}_{\mathcal{C}}} \leq 4$.
\end{definition}
\begin{definition}
	\label{def:wc}
	Two clusters in a level are said to be \textbf{well separated clusters}, if their respective boxes do not share boundary. Then, $\mathcal{W}_{\mathcal{C}}$ for a cluster $\mathcal{C}$ be a set containing clusters that are well separated from $\mathcal{C}$.
\end{definition}
For a cluster $\mathcal{C}$ corresponding to the node $i$ at level $l$, $\mathcal{N}^{(l)}_i$, we classify other clusters at the same level using the three disjoint sets, namely, the set of Edge sharers ($\mathcal{E}_\mathcal{C}$), the set of Vertex sharers ($\mathcal{V}_\mathcal{C}$) and the set of Well separated clusters ($\mathcal{W}_\mathcal{C}$). 
From the above definitions, the cluster $\mathcal{C}_0^{(0)}$ corresponding to the root node $\mathcal{N}^{(0)}_0$ is given as,
\begin{equation}
	{\mathcal{C}_0^{(0)}} = {\mathcal{C}}\cup\mathcal{E}_{\mathcal{C}} \cup \mathcal{V}_{\mathcal{C}} \cup \mathcal{W}_{\mathcal{C}} 
\end{equation}
Figure \ref{fig:clusRel} illustrates different clusters for a cluster (coloured red) at level 2. 

\begin{definition}
	\label{def:cc}
	\textbf{Clan set} $C_{\mathcal{C}}$ - For a cluster $\mathcal{C}$, clan set contains the $\text{siblings}(\mathcal{C})$ and the children of clusters in parents set of edge sharers (i.e., $\mathcal{E}_{\text{parent}(\mathcal{C})}$).
\end{definition}
The clan set $C_\mathcal{C}$ keeps track of the clusters that are unaccounted for in previous levels (i.e., the children of the parent's edge-sharing clusters and their siblings). In HODLR2D, we use the above defined four sets (i.e., $\mathcal{E}_\mathcal{C}$, $\mathcal{V}_\mathcal{C}$, $\mathcal{W}_\mathcal{C}$ and  $C_\mathcal{C}$) for each node in the hierarchical tree.
\begin{remark}
	\textbf{Admissibility criteria for HODLR2D} - Two clusters are admissible, if their respective boxes in the quadtree don't share an edge.
\end{remark} 
As with any hierarchical low-rank representations, the Interaction list of a cluster $\mathcal{C}$ ($I_\mathcal{C}$)  at a level has the clusters that are admissible (i.e., their interaction can be represented as low-rank) and that are not accounted for at previous level. 

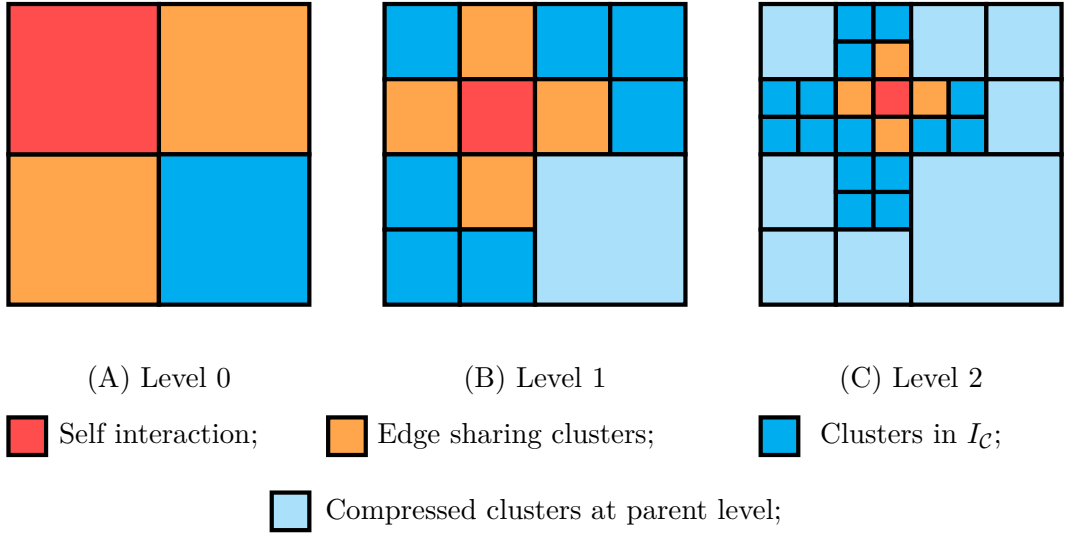
\begin{figure}[H]
	\centering
	\begin{tikzpicture}
		[%%%%%%%%%%%%%%%%%%%%%%%%%%%%%%
		box/.style={rectangle,draw=black, minimum size=1cm},
		]%%%%%%%%%%%%%%%%%%%%%%%%%%%%%%
		
		% Level 0
		\rect{orange!70}{0,0}{2,2}
		\rect{orange!70}{4,4}{2,2}
		\rect{red!70}{0,4}{2,2}
		\rect{cyan}{2,2}{4,0}
		
		\node at (2,-1) {(A) Level $0$};
		% Level 1
		
		\rect{cyan}{5,3}{6,4}
		\rect{orange!70}{6,3}{7,4}
		\rect{cyan}{7,3}{8,4}
		\rect{cyan}{8,3}{9,4}
		
		\rect{orange!70}{5,2}{6,3}
		\rect{red!70}{6,2}{7,3}
		\rect{orange!70}{7,2}{8,3}
		\rect{cyan}{8,2}{9,3}
		
		\rect{cyan}{5,1}{6,2}
		\rect{orange!70}{6,1}{7,2}
		\rect{cyan}{5,0}{6,1}
		\rect{cyan}{6,0}{7,1}
		
		\rect{cyan!30}{7,2}{9,0}
		
		\node at (7,-1) {(B) Level $1$};
		
		% Level 2
		
		\rect{cyan!30}{10,3}{11,4}
		
		\rect{cyan}{11,3}{11.5,3.5}
		\rect{cyan}{11,3.5}{11.5,4}
		\rect{orange!70}{11.5,3}{12,3.5}
		\rect{cyan}{11.5,3.5}{12,4}
		
		\rect{cyan!30}{12,3}{13,4}
		\rect{cyan!30}{13,3}{14,4}
		
		\rect{cyan}{10,2}{10.5,2.5}
		\rect{cyan}{10,2.5}{10.5,3}
		\rect{cyan}{10.5,2}{11,2.5}
		\rect{cyan}{10.5,2.5}{11,3}
		
		\rect{cyan}{11,2}{11.5,3.5}
		\rect{orange!70}{11,2.5}{11.5,3}
		\rect{orange!70}{11.5,2}{12,3.5}
		\rect{red!70}{11.5,2.5}{12,3}
		
		\rect{cyan}{12,2}{12.5,2.5}
		\rect{cyan}{12.5,2}{13,2.5}
		\rect{orange!70}{12,2.5}{12.5,3}
		\rect{cyan}{12.5,2.5}{13,3}
		
		\rect{cyan!30}{13,2}{14,3}
		
		\rect{cyan!30}{10,1}{11,2}
		
		\rect{cyan}{11,1}{11.5,1.5}
		\rect{cyan}{11,1.5}{11.5,2}
		\rect{cyan}{11.5,1}{12,1.5}
		\rect{cyan}{11.5,1.5}{12,2}
		
		\rect{cyan!30}{10,0}{11,1}
		\rect{cyan!30}{11,0}{12,1}
		
		\rect{cyan!30}{12,2}{14,0}
		
		\node at (12,-1) {(C) Level $2$};
		
		\rect{red!70}{0,-2}{0.5,-1.5};
		\node at (2,-1.75) {Self interaction;};
		
		\rect{orange!70}{4.25,-2}{4.75,-1.5};
		\node at (6.75,-1.75) {Edge sharing clusters;};
		
		\rect{cyan}{10,-2}{10.5,-1.5};
		\node at (12,-1.75) {Clusters in $I_\mathcal{C}$;};
		
		\rect{cyan!30}{3.5,-3}{4,-2.5};
		\node at (7.25,-2.75) {Compressed clusters at parent level;};
	\end{tikzpicture}
	\caption{Interaction relation of a cluster with other clusters for different levels.}
	\label{fig:three dom}
\end{figure}

\begin{definition}
	\label{def:ic}
	For a cluster $\mathcal{C}, \mathcal{I}_\mathcal{C}$ denotes a set defined by Equation~\eqref{eq:ilist}.
	\begin{equation}
		\label{eq:ilist}
		\mathcal{I}_\mathcal{C} = (C_\mathcal{C}\cap\mathcal{V}_\mathcal{C} )\cup (C_\mathcal{C} \cap\mathcal{W}_\mathcal{C}) 
	\end{equation}
	Then, the \textbf{Interaction list} ($I_\mathcal{C}$) of the cluster $\mathcal{C}$, is the list formed from the set $\mathcal{I}_\mathcal{C}$.
\end{definition} 
With above definitions, given a cluster $\mathcal{C}$ at level $l$ in the hierarchical tree, its interaction list is shown in Figure \ref{fig:three dom}.

\subsection{HODLR2D Algorithm}
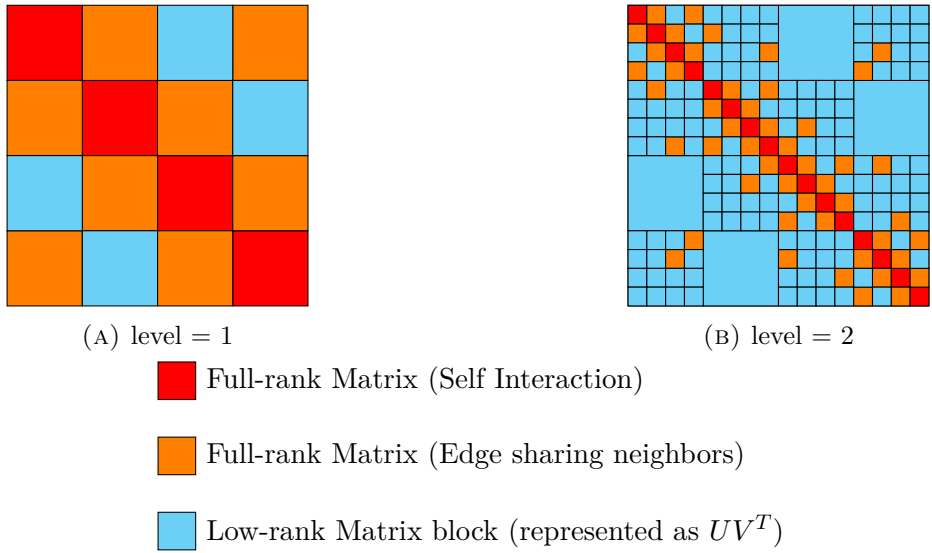
\begin{figure}[H]
	\begin{subfigure}[b]{0.5\textwidth}
		\centering
		\begin{tikzpicture}
			[%%%%%%%%%%%%%%%%%%%%%%%%%%%%%%
			box/.style={rectangle,draw=black, minimum size=0.5cm},
			]%%%%%%%%%%%%%%%%%%%%%%%%%%%%%%
			% DIAGONAL BLOCKS
			\draw[fill=red] (0,3) rectangle +(1,1);
			\draw[fill=red] (1,2) rectangle +(1,1);
			\draw[fill=red] (2,1) rectangle +(1,1);
			\draw[fill=red] (3,0) rectangle +(1,1);
			%EDGE SHARING BLOCKS
			\draw[fill=orange] (2,0) rectangle +(1,1);
			\draw[fill=orange] (0,2) rectangle +(1,1);
			
			\draw[fill=orange] (0,0) rectangle +(1,1);
			\draw[fill=orange] (1,1) rectangle +(1,1);
			\draw[fill=orange] (2,2) rectangle +(1,1);
			\draw[fill=orange] (3,3) rectangle +(1,1);
			
			\draw[fill=orange] (1,3) rectangle +(1,1);
			\draw[fill=orange] (3,1) rectangle +(1,1);
			% VERTEX SHARING BLOCKS
			\draw[fill=cyan!50] (2,3) rectangle +(1,1);
			\draw[fill=cyan!50] (3,2) rectangle +(1,1);
			\draw[fill=cyan!50] (1,0) rectangle +(1,1);
			\draw[fill=cyan!50] (0,1) rectangle +(1,1);
			
		\end{tikzpicture}
		\caption{level = 1}
	\end{subfigure}%
	\hfill
	\begin{subfigure}[b]{0.5\textwidth}
		\centering
		\begin{tikzpicture}[scale=0.4]
			[%%%%%%%%%%%%%%%%%%%%%%%%%%%%%%
			box/.style={rectangle,draw=black, minimum size=0.5cm},
			]%%%%%%%%%%%%%%%%%%%%%%%%%%%%%%
			\draw (0,0) rectangle +(10,10);
			\draw[draw=black,fill=cyan!50] (0,2.5)rectangle +(2.5,2.5);
			\draw[draw=black,fill=cyan!50] (2.5,0)rectangle +(2.5,2.5);
			\draw[draw=black,fill=cyan!50] (5,7.5)rectangle +(2.5,2.5);
			\draw[draw=black,fill=cyan!50] (7.5,5)rectangle +(2.5,2.5);
			\draw[draw=black,fill=red] (0,9.375)rectangle +(0.625,0.625);
			\draw[draw=black,fill=orange](0,8.75) rectangle +(0.625,0.625);
			\draw[draw=black,fill=orange](0,7.5) rectangle +(0.625,0.625);
			\draw[draw=black,fill=cyan!50] (0,8.125)rectangle +(0.625,0.625);
			\draw[draw=black,fill=cyan!50] (0,6.875)rectangle +(0.625,0.625);
			\draw[draw=black,fill=cyan!50] (0,6.25)rectangle +(0.625,0.625);
			\draw[draw=black,fill=cyan!50] (0,5.625)rectangle +(0.625,0.625);
			\draw[draw=black,fill=cyan!50] (0,5)rectangle +(0.625,0.625);
			\draw[draw=black,fill=cyan!50] (0,1.875)rectangle +(0.625,0.625);
			\draw[draw=black,fill=cyan!50] (0,1.25)rectangle +(0.625,0.625);
			\draw[draw=black,fill=cyan!50] (0,0.625)rectangle +(0.625,0.625);
			\draw[draw=black,fill=cyan!50] (0,0)rectangle +(0.625,0.625);
			\draw[draw=black,fill=red] (0.625,8.75)rectangle +(0.625,0.625);
			\draw[draw=black,fill=orange](0.625,9.375) rectangle +(0.625,0.625);
			\draw[draw=black,fill=orange](0.625,6.875) rectangle +(0.625,0.625);
			\draw[draw=black,fill=orange](0.625,8.125) rectangle +(0.625,0.625);
			\draw[draw=black,fill=cyan!50] (0.625,7.5)rectangle +(0.625,0.625);
			\draw[draw=black,fill=cyan!50] (0.625,6.25)rectangle +(0.625,0.625);
			\draw[draw=black,fill=cyan!50] (0.625,5.625)rectangle +(0.625,0.625);
			\draw[draw=black,fill=cyan!50] (0.625,5)rectangle +(0.625,0.625);
			\draw[draw=black,fill=cyan!50] (0.625,1.875)rectangle +(0.625,0.625);
			\draw[draw=black,fill=cyan!50] (0.625,1.25)rectangle +(0.625,0.625);
			\draw[draw=black,fill=cyan!50] (0.625,0.625)rectangle +(0.625,0.625);
			\draw[draw=black,fill=cyan!50] (0.625,0)rectangle +(0.625,0.625);
			\draw[draw=black,fill=red] (1.25,8.125)rectangle +(0.625,0.625);
			\draw[draw=black,fill=orange](1.25,8.75) rectangle +(0.625,0.625);
			\draw[draw=black,fill=orange](1.25,7.5) rectangle +(0.625,0.625);
			\draw[draw=black,fill=orange](1.25,5) rectangle +(0.625,0.625);
			\draw[draw=black,fill=orange](1.25,1.25) rectangle +(0.625,0.625);
			\draw[draw=black,fill=cyan!50] (1.25,9.375)rectangle +(0.625,0.625);
			\draw[draw=black,fill=cyan!50] (1.25,6.875)rectangle +(0.625,0.625);
			\draw[draw=black,fill=cyan!50] (1.25,6.25)rectangle +(0.625,0.625);
			\draw[draw=black,fill=cyan!50] (1.25,5.625)rectangle +(0.625,0.625);
			\draw[draw=black,fill=cyan!50] (1.25,1.875)rectangle +(0.625,0.625);
			\draw[draw=black,fill=cyan!50] (1.25,0.625)rectangle +(0.625,0.625);
			\draw[draw=black,fill=cyan!50] (1.25,0)rectangle +(0.625,0.625);
			\draw[draw=black,fill=red] (1.875,7.5)rectangle +(0.625,0.625);
			\draw[draw=black,fill=orange](1.875,9.375) rectangle +(0.625,0.625);
			\draw[draw=black,fill=orange](1.875,8.125) rectangle +(0.625,0.625);
			\draw[draw=black,fill=orange](1.875,1.875) rectangle +(0.625,0.625);
			\draw[draw=black,fill=cyan!50] (1.875,8.75)rectangle +(0.625,0.625);
			\draw[draw=black,fill=cyan!50] (1.875,6.875)rectangle +(0.625,0.625);
			\draw[draw=black,fill=cyan!50] (1.875,6.25)rectangle +(0.625,0.625);
			\draw[draw=black,fill=cyan!50] (1.875,5.625)rectangle +(0.625,0.625);
			\draw[draw=black,fill=cyan!50] (1.875,5)rectangle +(0.625,0.625);
			\draw[draw=black,fill=cyan!50] (1.875,1.25)rectangle +(0.625,0.625);
			\draw[draw=black,fill=cyan!50] (1.875,0.625)rectangle +(0.625,0.625);
			\draw[draw=black,fill=cyan!50] (1.875,0)rectangle +(0.625,0.625);
			\draw[draw=black,fill=red] (2.5,6.875)rectangle +(0.625,0.625);
			\draw[draw=black,fill=orange](2.5,8.75) rectangle +(0.625,0.625);
			\draw[draw=black,fill=orange](2.5,6.25) rectangle +(0.625,0.625);
			\draw[draw=black,fill=orange](2.5,5) rectangle +(0.625,0.625);
			\draw[draw=black,fill=cyan!50] (2.5,5.625)rectangle +(0.625,0.625);
			\draw[draw=black,fill=cyan!50] (2.5,9.375)rectangle +(0.625,0.625);
			\draw[draw=black,fill=cyan!50] (2.5,8.125)rectangle +(0.625,0.625);
			\draw[draw=black,fill=cyan!50] (2.5,7.5)rectangle +(0.625,0.625);
			\draw[draw=black,fill=cyan!50] (2.5,4.375)rectangle +(0.625,0.625);
			\draw[draw=black,fill=cyan!50] (2.5,3.75)rectangle +(0.625,0.625);
			\draw[draw=black,fill=cyan!50] (2.5,3.125)rectangle +(0.625,0.625);
			\draw[draw=black,fill=cyan!50] (2.5,2.5)rectangle +(0.625,0.625);
			\draw[draw=black,fill=red] (3.125,6.25)rectangle +(0.625,0.625);
			\draw[draw=black,fill=orange](3.125,6.875) rectangle +(0.625,0.625);
			\draw[draw=black,fill=orange](3.125,5.625) rectangle +(0.625,0.625);
			\draw[draw=black,fill=cyan!50] (3.125,5)rectangle +(0.625,0.625);
			\draw[draw=black,fill=cyan!50] (3.125,9.375)rectangle +(0.625,0.625);
			\draw[draw=black,fill=cyan!50] (3.125,8.75)rectangle +(0.625,0.625);
			\draw[draw=black,fill=cyan!50] (3.125,8.125)rectangle +(0.625,0.625);
			\draw[draw=black,fill=cyan!50] (3.125,7.5)rectangle +(0.625,0.625);
			\draw[draw=black,fill=cyan!50] (3.125,4.375)rectangle +(0.625,0.625);
			\draw[draw=black,fill=cyan!50] (3.125,3.75)rectangle +(0.625,0.625);
			\draw[draw=black,fill=cyan!50] (3.125,3.125)rectangle +(0.625,0.625);
			\draw[draw=black,fill=cyan!50] (3.125,2.5)rectangle +(0.625,0.625);
			\draw[draw=black,fill=red] (3.75,5.625)rectangle +(0.625,0.625);
			\draw[draw=black,fill=orange](3.75,6.25) rectangle +(0.625,0.625);
			\draw[draw=black,fill=orange](3.75,5) rectangle +(0.625,0.625);
			\draw[draw=black,fill=orange](3.75,3.75) rectangle +(0.625,0.625);
			\draw[draw=black,fill=cyan!50] (3.75,6.875)rectangle +(0.625,0.625);
			\draw[draw=black,fill=cyan!50] (3.75,9.375)rectangle +(0.625,0.625);
			\draw[draw=black,fill=cyan!50] (3.75,8.75)rectangle +(0.625,0.625);
			\draw[draw=black,fill=cyan!50] (3.75,8.125)rectangle +(0.625,0.625);
			\draw[draw=black,fill=cyan!50] (3.75,7.5)rectangle +(0.625,0.625);
			\draw[draw=black,fill=cyan!50] (3.75,4.375)rectangle +(0.625,0.625);
			\draw[draw=black,fill=cyan!50] (3.75,3.125)rectangle +(0.625,0.625);
			\draw[draw=black,fill=cyan!50] (3.75,2.5)rectangle +(0.625,0.625);
			\draw[draw=black,fill=red] (4.375,5)rectangle +(0.625,0.625);
			\draw[draw=black,fill=orange](4.375,6.875) rectangle +(0.625,0.625);
			\draw[draw=black,fill=orange](4.375,8.125) rectangle +(0.625,0.625);
			\draw[draw=black,fill=orange](4.375,5.625) rectangle +(0.625,0.625);
			\draw[draw=black,fill=orange](4.375,4.375) rectangle +(0.625,0.625);
			\draw[draw=black,fill=cyan!50] (4.375,6.25)rectangle +(0.625,0.625);
			\draw[draw=black,fill=cyan!50] (4.375,9.375)rectangle +(0.625,0.625);
			\draw[draw=black,fill=cyan!50] (4.375,8.75)rectangle +(0.625,0.625);
			\draw[draw=black,fill=cyan!50] (4.375,7.5)rectangle +(0.625,0.625);
			\draw[draw=black,fill=cyan!50] (4.375,3.75)rectangle +(0.625,0.625);
			\draw[draw=black,fill=cyan!50] (4.375,3.125)rectangle +(0.625,0.625);
			\draw[draw=black,fill=cyan!50] (4.375,2.5)rectangle +(0.625,0.625);
			\draw[draw=black,fill=red] (5,4.375)rectangle +(0.625,0.625);
			\draw[draw=black,fill=orange](5,5) rectangle +(0.625,0.625);
			\draw[draw=black,fill=orange](5,1.25) rectangle +(0.625,0.625);
			\draw[draw=black,fill=orange](5,3.75) rectangle +(0.625,0.625);
			\draw[draw=black,fill=orange](5,2.5) rectangle +(0.625,0.625);
			\draw[draw=black,fill=cyan!50] (5,3.125)rectangle +(0.625,0.625);
			\draw[draw=black,fill=cyan!50] (5,6.875)rectangle +(0.625,0.625);
			\draw[draw=black,fill=cyan!50] (5,6.25)rectangle +(0.625,0.625);
			\draw[draw=black,fill=cyan!50] (5,5.625)rectangle +(0.625,0.625);
			\draw[draw=black,fill=cyan!50] (5,1.875)rectangle +(0.625,0.625);
			\draw[draw=black,fill=cyan!50] (5,0.625)rectangle +(0.625,0.625);
			\draw[draw=black,fill=cyan!50] (5,0)rectangle +(0.625,0.625);
			\draw[draw=black,fill=red] (5.625,3.75)rectangle +(0.625,0.625);
			\draw[draw=black,fill=orange](5.625,5.625) rectangle +(0.625,0.625);
			\draw[draw=black,fill=orange](5.625,4.375) rectangle +(0.625,0.625);
			\draw[draw=black,fill=orange](5.625,3.125) rectangle +(0.625,0.625);
			\draw[draw=black,fill=cyan!50] (5.625,2.5)rectangle +(0.625,0.625);
			\draw[draw=black,fill=cyan!50] (5.625,6.875)rectangle +(0.625,0.625);
			\draw[draw=black,fill=cyan!50] (5.625,6.25)rectangle +(0.625,0.625);
			\draw[draw=black,fill=cyan!50] (5.625,5)rectangle +(0.625,0.625);
			\draw[draw=black,fill=cyan!50] (5.625,1.875)rectangle +(0.625,0.625);
			\draw[draw=black,fill=cyan!50] (5.625,1.25)rectangle +(0.625,0.625);
			\draw[draw=black,fill=cyan!50] (5.625,0.625)rectangle +(0.625,0.625);
			\draw[draw=black,fill=cyan!50] (5.625,0)rectangle +(0.625,0.625);
			\draw[draw=black,fill=red] (6.25,3.125)rectangle +(0.625,0.625);
			\draw[draw=black,fill=orange](6.25,3.75) rectangle +(0.625,0.625);
			\draw[draw=black,fill=orange](6.25,2.5) rectangle +(0.625,0.625);
			\draw[draw=black,fill=cyan!50] (6.25,4.375)rectangle +(0.625,0.625);
			\draw[draw=black,fill=cyan!50] (6.25,6.875)rectangle +(0.625,0.625);
			\draw[draw=black,fill=cyan!50] (6.25,6.25)rectangle +(0.625,0.625);
			\draw[draw=black,fill=cyan!50] (6.25,5.625)rectangle +(0.625,0.625);
			\draw[draw=black,fill=cyan!50] (6.25,5)rectangle +(0.625,0.625);
			\draw[draw=black,fill=cyan!50] (6.25,1.875)rectangle +(0.625,0.625);
			\draw[draw=black,fill=cyan!50] (6.25,1.25)rectangle +(0.625,0.625);
			\draw[draw=black,fill=cyan!50] (6.25,0.625)rectangle +(0.625,0.625);
			\draw[draw=black,fill=cyan!50] (6.25,0)rectangle +(0.625,0.625);
			\draw[draw=black,fill=red] (6.875,2.5)rectangle +(0.625,0.625);
			\draw[draw=black,fill=orange](6.875,4.375) rectangle +(0.625,0.625);
			\draw[draw=black,fill=orange](6.875,0.625) rectangle +(0.625,0.625);
			\draw[draw=black,fill=orange](6.875,3.125) rectangle +(0.625,0.625);
			\draw[draw=black,fill=cyan!50] (6.875,3.75)rectangle +(0.625,0.625);
			\draw[draw=black,fill=cyan!50] (6.875,6.875)rectangle +(0.625,0.625);
			\draw[draw=black,fill=cyan!50] (6.875,6.25)rectangle +(0.625,0.625);
			\draw[draw=black,fill=cyan!50] (6.875,5.625)rectangle +(0.625,0.625);
			\draw[draw=black,fill=cyan!50] (6.875,5)rectangle +(0.625,0.625);
			\draw[draw=black,fill=cyan!50] (6.875,1.875)rectangle +(0.625,0.625);
			\draw[draw=black,fill=cyan!50] (6.875,1.25)rectangle +(0.625,0.625);
			\draw[draw=black,fill=cyan!50] (6.875,0)rectangle +(0.625,0.625);
			\draw[draw=black,fill=red] (7.5,1.875)rectangle +(0.625,0.625);
			\draw[draw=black,fill=orange](7.5,7.5) rectangle +(0.625,0.625);
			\draw[draw=black,fill=orange](7.5,1.25) rectangle +(0.625,0.625);
			\draw[draw=black,fill=orange](7.5,0) rectangle +(0.625,0.625);
			\draw[draw=black,fill=cyan!50] (7.5,0.625)rectangle +(0.625,0.625);
			\draw[draw=black,fill=cyan!50] (7.5,9.375)rectangle +(0.625,0.625);
			\draw[draw=black,fill=cyan!50] (7.5,8.75)rectangle +(0.625,0.625);
			\draw[draw=black,fill=cyan!50] (7.5,8.125)rectangle +(0.625,0.625);
			\draw[draw=black,fill=cyan!50] (7.5,4.375)rectangle +(0.625,0.625);
			\draw[draw=black,fill=cyan!50] (7.5,3.75)rectangle +(0.625,0.625);
			\draw[draw=black,fill=cyan!50] (7.5,3.125)rectangle +(0.625,0.625);
			\draw[draw=black,fill=cyan!50] (7.5,2.5)rectangle +(0.625,0.625);
			\draw[draw=black,fill=red] (8.125,1.25)rectangle +(0.625,0.625);
			\draw[draw=black,fill=orange](8.125,8.125) rectangle +(0.625,0.625);
			\draw[draw=black,fill=orange](8.125,1.875) rectangle +(0.625,0.625);
			\draw[draw=black,fill=orange](8.125,4.375) rectangle +(0.625,0.625);
			\draw[draw=black,fill=orange](8.125,0.625) rectangle +(0.625,0.625);
			\draw[draw=black,fill=cyan!50] (8.125,0)rectangle +(0.625,0.625);
			\draw[draw=black,fill=cyan!50] (8.125,9.375)rectangle +(0.625,0.625);
			\draw[draw=black,fill=cyan!50] (8.125,8.75)rectangle +(0.625,0.625);
			\draw[draw=black,fill=cyan!50] (8.125,7.5)rectangle +(0.625,0.625);
			\draw[draw=black,fill=cyan!50] (8.125,3.75)rectangle +(0.625,0.625);
			\draw[draw=black,fill=cyan!50] (8.125,3.125)rectangle +(0.625,0.625);
			\draw[draw=black,fill=cyan!50] (8.125,2.5)rectangle +(0.625,0.625);
			\draw[draw=black,fill=red] (8.75,0.625)rectangle +(0.625,0.625);
			\draw[draw=black,fill=orange](8.75,1.25) rectangle +(0.625,0.625);
			\draw[draw=black,fill=orange](8.75,0) rectangle +(0.625,0.625);
			\draw[draw=black,fill=orange](8.75,2.5) rectangle +(0.625,0.625);
			\draw[draw=black,fill=cyan!50] (8.75,1.875)rectangle +(0.625,0.625);
			\draw[draw=black,fill=cyan!50] (8.75,9.375)rectangle +(0.625,0.625);
			\draw[draw=black,fill=cyan!50] (8.75,8.75)rectangle +(0.625,0.625);
			\draw[draw=black,fill=cyan!50] (8.75,8.125)rectangle +(0.625,0.625);
			\draw[draw=black,fill=cyan!50] (8.75,7.5)rectangle +(0.625,0.625);
			\draw[draw=black,fill=cyan!50] (8.75,4.375)rectangle +(0.625,0.625);
			\draw[draw=black,fill=cyan!50] (8.75,3.75)rectangle +(0.625,0.625);
			\draw[draw=black,fill=cyan!50] (8.75,3.125)rectangle +(0.625,0.625);
			\draw[draw=black,fill=red] (9.375,0)rectangle +(0.625,0.625);
			\draw[draw=black,fill=orange](9.375,1.875) rectangle +(0.625,0.625);
			\draw[draw=black,fill=orange](9.375,0.625) rectangle +(0.625,0.625);
			\draw[draw=black,fill=cyan!50] (9.375,1.25)rectangle +(0.625,0.625);
			\draw[draw=black,fill=cyan!50] (9.375,9.375)rectangle +(0.625,0.625);
			\draw[draw=black,fill=cyan!50] (9.375,8.75)rectangle +(0.625,0.625);
			\draw[draw=black,fill=cyan!50] (9.375,8.125)rectangle +(0.625,0.625);
			\draw[draw=black,fill=cyan!50] (9.375,7.5)rectangle +(0.625,0.625);
			\draw[draw=black,fill=cyan!50] (9.375,4.375)rectangle +(0.625,0.625);
			\draw[draw=black,fill=cyan!50] (9.375,3.75)rectangle +(0.625,0.625);
			\draw[draw=black,fill=cyan!50] (9.375,3.125)rectangle +(0.625,0.625);
			\draw[draw=black,fill=cyan!50] (9.375,2.5)rectangle +(0.625,0.625);
		\end{tikzpicture}
		\caption{level = 2}
	\end{subfigure}%
	\hfill
	\begin{subfigure}[b]{0.5\textwidth}
		\centering
		\begin{tikzpicture}
			[%%%%%%%%%%%%%%%%%%%%%%%%%%%%%%
			box/.style={rectangle,draw=black, minimum size=0.5cm},
			]%%%%%%%%%%%%%%%%%%%%%%%%%%%%%%
			
			\node[box,fill=red,label=right:Full-rank Matrix (Self Interaction),anchor=west] at (0,2){};
			\node[box,fill=orange,label=right:Full-rank Matrix (Edge sharing neighbors),anchor=west] at (0,1){};
			\node[box,fill=cyan!50,label=right:Low-rank Matrix block (represented as $UV^T$),anchor=west] at (0,0){};
		\end{tikzpicture}
		\caption*{}
	\end{subfigure}%
	\caption{HODLR2D matrix at different levels.}
	\label{fig:hodlr2d_mat}
\end{figure}

Algorithm \ref{alg:cap} provides the initialization routine of the HODLR2D based on quadtree. We start with the subdivision of the box corresponding to each the cluster into four geometrically disjoint boxes and number the resulting smaller boxes as in Figure \ref{sdomain01}. The Box $B\subset \Rb^2$ is subdivided using a $\kappa$ level balanced quadtree such that all boxes at the leaf level have at most $N_{\max}$ particles. Then for all the nodes in the Hierarchical tree, we form the set of Edge sharers ($\mathcal{E}_\mathcal{C}$) and the Interaction list ($I_\mathcal{C}$). At all levels, we compress and represent the matrices corresponding to the elements in the $I_\mathcal{C}$ using ACA. The resulting low-rank approximation of the matrix at different levels using HODLR2D is as in Figure \ref{fig:hodlr2d_mat}.

\begin{remark}
	In HODLR2D, the edge sharing clusters at the leaf level are represented as full-rank matrix blocks even though they are numerically low-rank.
\end{remark} 

%% HODLR2D MATRIX at Different Levels

Upon constructing the HODLR2D representation using Algorithm~\ref{alg:cap}, we accelerate the matrix-vector product using Algorithm~\ref{alg:matvec}. In Algorithm \ref{alg:matvec}, lines $4$ and $7$ represent a dense matrix-vector operation corresponding to self-interaction and interaction between edge-sharing clusters at the leaf level, respectively. Moreover, line $15$ performs a low-rank matrix-vector product for all nodes in the quadtree. For a node in the hierarchical tree of HODLR2D at any level, the maximum size of the interaction list is $15$. At the leaf level, a node can have a maximum of five dense matrices (including the self-interaction). Comparing this with the $\mathcal{H}$-matrix described in Section~\ref{sec2}, the maximum size of the interaction list is $27$, and at the leaf level, a node can have a maximum of nine dense matrices. Furthermore, another difference is that the low-rank approximation for HODLR2D starts at nodes in level $1$, whereas for $\mathcal{H}$-matrix, it starts at level $2$.

\begin{algorithm}[H]
	\caption{HODLR2D}\label{alg:cap}
	\begin{algorithmic}[1]
		\Procedure{InitializeHODLR2D}{$N_{\max}$,$\epsilon$}
		\State{} \Comment{$N_{\max}$ is the maximum number of particles at leaf level; $\epsilon$ is the tolerance for ACA}
		\State {Consider a balanced quadtree with $\kappa$ levels such that each node at the leaf level contains not more than $N_{\max}$ particles.
		}
		\State {For each nodes in quadtree from root to leaf compute Interaction list $I_\mathcal{C}$ and set of Edge sharers $\mathcal{E}_\mathcal{C}$ using definitions~\ref{def:ic} and ~\ref{def:ec} respectively.}
		\State {For each nodes at leaf level $\kappa$, compute the dense matrix corresponding to the self interaction and the dense matrices corresponding to edge sharing list.}
		\State {For each level $l \in\{\kappa,\kappa-1,\ldots ,1\}$, and for all nodes in a particular level, compute the low rank approximation $UV^T$ using ACA with prescribed tolerance $\epsilon$ corresponding to the clusters in Interaction list $I_\mathcal{C}$.}
		\EndProcedure
	\end{algorithmic}
\end{algorithm}

\begin{algorithm}[H]
	\caption{HODLR2D Matrix Vector Product $K\psi = b$}\label{alg:matvec}
	\begin{algorithmic}[1]
		\Procedure{MatVec}{$\psi$}
		\For{\texttt{i=1:$4^\kappa$}} \Comment{Full-rank Mat-Vec product}
		\State $X\gets \text{Index set of } \mathcal{C}_{i}^{(\kappa)}$
		\State $b(X) = b(X) + K(X,X) \times \psi(X)$
		\For{\texttt{j in $\mathcal{E}_\mathcal{C}$}}
		\State $Y\gets \text{Index set of } \mathcal{C}_{j}^{(\kappa)}$
		\State $b(X) = b(X) + K(X,Y) \times \psi(Y)$
		\EndFor
		\EndFor
		
		\For{\texttt{$l=1:\kappa$}} \Comment{Low-rank Mat-Vec product}
		\For{\texttt{i=1:$4^\kappa$}}
		\State $X\gets \text{Index set of } \mathcal{C}_{i}^{(l)}$
		\For{\texttt{j in $I_\mathcal{C}$}}
		\State $Y\gets \text{Index set of } \mathcal{C}_{j}^{(l)}$
		\State $b(X) = b(X) + U_{ij}^{(l)}\times \bkt{{V_{ij}^{(l)}}^T\times \psi(Y)}$
		\EndFor
		\EndFor
		\EndFor
		\State \textbf{return} $b$
		\EndProcedure
	\end{algorithmic}
\end{algorithm}

The total storage cost, initialization time and computational cost to perform matrix-vector products scale as $\mathcal{O}\bkt{pN \kappa}$, where $p$ is the maximum rank of vertex sharing blocks and $\kappa$ is the number of levels in the balanced quadtree. From Theorem~\ref{tm:farfield}, we have $p = \mathcal{O}\bkt{\log \bkt{N} \log\bkt{\log\bkt{N}/\epsilon}}$ and $\kappa \in \mathcal{O}\bkt{\log \bkt{N}}$ for a uniform distribution of particles. Hence, the total storage cost and computational cost to perform matrix-vector products scale as $\mathcal{O}\bkt{N \log^2\bkt{N}\log\bkt{\log\bkt{N}/\epsilon}}$. Note that the algorithm can also be adapted to a non-balanced or adaptive quadtree with little modifications. We work on a balanced quadtree for pedagogical reasons.
	\section{Numerical Experiments}
\label{sec5}
This section demonstrates the performance of the HODLR2D format and compares it with HODLR and $\mathcal{H}$-matrix formats. To being with, we compare the scaling of matrix-vector products of HODLR2D matrix with HODLR and $\mathcal{H}$-matrix. We then accelerate GMRES based iterative solver for dense linear systems using HODLR2D, HOLDR and $\mathcal{H}$-matrices. The first application involves a dense linear system arising out of radial basis function interpolation. The next application involves a dense linear system arising out the discretization of an integral equation. In both examples, we use HODLR, HODLR2D and $\mathcal{H}$-matrix formats to approximate the underlying dense matrix and use it to accelerate the matrix-vector product in GMRES. We use ACA to perform the low-rank decomposition in all three hierarchical representations. The considered example problems are as generic as possible. Throughout our experiments, the following parameters remain unchanged.
\begin{itemize}
	\item Tolerance for computing low-rank approximation through ACA to $10^{-12}$
	\item Stopping criteria for GMRES is residual less than $10^{-10}$ with restart after each iterations.
	\item Number of levels in the hierarchical tree is decided such that the number of charges at leaf level does not exceed 500.
\end{itemize}
We tabulate the results of numerical experiments with the following notations.
\begin{table}[H]
	\resizebox{\textwidth}{!}{%
		\begin{tabular}{|l|l|}
			\hline
			$N$            & Number of Unknowns/Degrees of 
			Freedom\\ \hline
			$N_{\max}$  & Number of particles at the leaf level in the tree; kept as $500$ throughout the article.\\
			\hline
			$T_I$        & Time to initialize the hierarchical tree (in seconds)\\
			\hline
			$T_G$        & Time taken by GMRES to converge (in seconds)\\
			\hline
			CR           & Compression Ratio, which denotes the ratio of number of FP64 values in the hierarchical tree to $N^2$\\
			\hline
			$r_m$        & Maximum rank across the hierarchical tree\\
			\hline
			$\epsilon_r$ & Relative error in the solution \\
			\hline
			$\epsilon$ & Tolerance for the ACA; kept as $10^{-12}$ throughout the article.\\
			\hline
		\end{tabular}%
	}
	\label{tab:nexp}
\end{table}
\subsection{HODLR2D matrix-vector product}
\label{sec_HODLR2D_matvec}
We first demonstrate the space and computational performance of HODLR2D in comparison with $\mathcal{H}$-matrix and HODLR using the kernel defined in Equation~\eqref{eq:1byr}. Let us now consider $N$ points $\{x_i\}_{i=0}^N$ in the Box $B \subset [-1,1]^2$ ($\sqrt{N} \times \sqrt{N}$ Chebyshev grid on $[-1,1]^2$). The matrix $K\in \mathbb{R}^{N\times N}$ due to the interaction of points $x_i$ defined by Equation~\eqref{eq:1byr}.
\begin{equation}	\label{eq:1byr}
	K_{ij} =
	\begin{cases}
		\frac{1}{r_{ij}}       & \quad i \neq j\\
		0  & \quad i = j
	\end{cases}
\end{equation}
$$r_{ij} = ||x_i - x_j||_2, \text{where, } x_i, x_j \in B$$
Table~\ref{tab:1byr_init} shows the maximum rank across the hierarchical tree, the memory required and the initialization time of HODLR2D, $\mathcal{H}$-matrix and HODLR for different system size $N$. It is evident from the Table~\ref{tab:1byr_init}, that the maximum rank, memory required and the initialization time for HODLR2D is way lesser than HODLR. Also, the memroy required and initialization time for HODLR2D is comparable with $\mathcal{H}$-matrix.

To demonstrate the computational performance of HODLR2D, $\mathcal{H}$-matrix and HODLR, we use ten different input vectors $\psi \in \mathbb{R}^{N\times 1}$. We obtain the right-hand side vector $b \in \mathbb{R}^{N\times 1}$ using explicit matrix-vector product $K\psi = b$. Let $\hat{K}$ represent the hierarchical low-rank representation of $K$ and let the computed right hand side using the algorithm be $\hat{b} = \hat{K}\psi$. Table~\ref{tab:1byr_mat} compares the computational performance of HODLR2D, $\mathcal{H}$-matrix and HODLR using the average time taken for matrix-vector product on the ten different pair $\psi \in \Rb^{N \times 1}$ and the error is the maximum relative error $\epsilon_r$, i.e., $\epsilon_r = \displaystyle \max_{i=1}^n\dfrac{\abs{\hat{b}-b}_i}{ \abs{b}_i}$. The error reported in Table~\ref{tab:1byr_mat} is the maximum error over the ten right hand sides. All the computations in this section were performed using only a single core of Intel Xeon Gold 6248, 20-core, 2.5 GHz processor with memory of 192GB. 
\begin{table}[H]
	\centering
	\resizebox{\textwidth}{!}{%
		\begin{tabular}{lccccccccc}
			\hline
			\multirow{2}{*}{N} & \multicolumn{3}{c}{$r_m$}                                         & \multicolumn{3}{c}{Memory used (in GB)}                           & \multicolumn{3}{c}{$T_I$}                                                     \\ \cline{2-10} 
			& \textbf{HODLR2D} & \textbf{$\mathcal{H}$-matrix} & \textbf{HODLR} & \textbf{HODLR2D} & \textbf{$\mathcal{H}$-matrix} & \textbf{HODLR} & \textbf{HODLR2D} & \textbf{$\mathcal{H}$-matrix} & \textbf{HODLR} \\ \hline
			10000              & 113              & 56                            & 515            & 0.23             & 0.26                          & 0.25           & 2.4628           & 3.17731                       & 4.77726        \\
			22500              & 127              & 60                            & 797            & 0.69             & 0.78                          & 0.88           & 8.02592          & 9.71345                       & 27.9686        \\
			40000              & 138              & 59                            & 1051           & 1.46             & 1.65                          & 2.12           & 19.8064          & 23.3388                       & 88.4981        \\
			62500              & 145              & 63                            & 1308           & 2.53             & 2.87                          & 4.14           & 31.5305          & 35.9836                       & 209.583        \\
			90000              & 159              & 61                            & 1499           & 3.99             & 4.53                          & 7.26           & 57.6382          & 76.5599                       & 463.171        \\
			160000             & 165              & 62                            & 2028           & 7.96             & 9.04                          & 17.06          & 125.052          & 168.782                       & 1404.71        \\
			250000             & 180              & 62                            & 2596           & 13.52            & 15.31                         & 33.57          & 186.154          & 222.405                       & 3501.93        \\ \hline
		\end{tabular}
	}
	\caption{Space complexity of HODLR2D, $\mathcal{H}$-matrix and HODLR for matrix $K$ whose elements are defined by Equation~\eqref{eq:1byr}}
	\label{tab:1byr_init}
\end{table}

\begin{table}[H]
	\centering
	%\resizebox{\textwidth}{!}{%
	\begin{tabular}{@{}cccccccccc@{}}
		\toprule
		\multirow{2}{*}{N} & \multicolumn{3}{c}{Matrix-Vector product time (in s)} & \multicolumn{3}{c}{$\epsilon_r$} \\ \cmidrule(lr){2-4}\cmidrule(lr){5-7}
		& HODLR2D & $\mathcal{H}$-matrix  & HODLR & HODLR2D & $\mathcal{H}$-matrix  & HODLR  \\ \hline
		10000              & 0.0284  & 0.0307 & 0.0211               & 1.5$\times10^{-13}$   & 1$\times10^{-13}$     & 3.46$\times10^{-12}$  \\
		22500              & 0.0903  & 0.1015 & 0.0723               & 1.75$\times10^{-12}$  & 8.7$\times10^{-13}$   & 1.181$\times10^{-11}$ \\
		40000              & 0.1979  & 0.228  & 0.1744               & 7.23$\times10^{-12}$  & 6.49$\times10^{-12}$  & 2.232$\times10^{-11}$ \\
		62500              & 0.3208  & 0.374  & 0.3378               & 5.16$\times10^{-12}$  & 5.01$\times10^{-12}$  & 1.751$\times10^{-10}$ \\
		90000              & 0.557   & 0.6824 & 0.6264               & 3.88$\times10^{-12}$  & 3.02$\times10^{-12}$  & 9.09$\times10^{-10}$  \\
		160000             & 1.2835  & 1.6359 & 1.5465               & 2.049$\times10^{-11}$ & 2.027$\times10^{-11}$ & 4.03$\times10^{-10}$  \\
		250000             & 1.9681  & 2.4578 & 3.055                & 2.588$\times10^{-11}$ & 2.132$\times10^{-11}$ & 2.037$\times10^{-9}$  \\ \hline
	\end{tabular}
	%}
	\caption{Computational performance of HODLR2D, $\mathcal{H}$-matrix and HODLR for matrix $K$ whose elements are defined by Equation~\eqref{eq:1byr}}
	\label{tab:1byr_mat}
\end{table}

\begin{figure}[H]
	\centering
	\begin{subfigure}{.5\textwidth}
		\centering
		\includegraphics[width=0.8\linewidth]{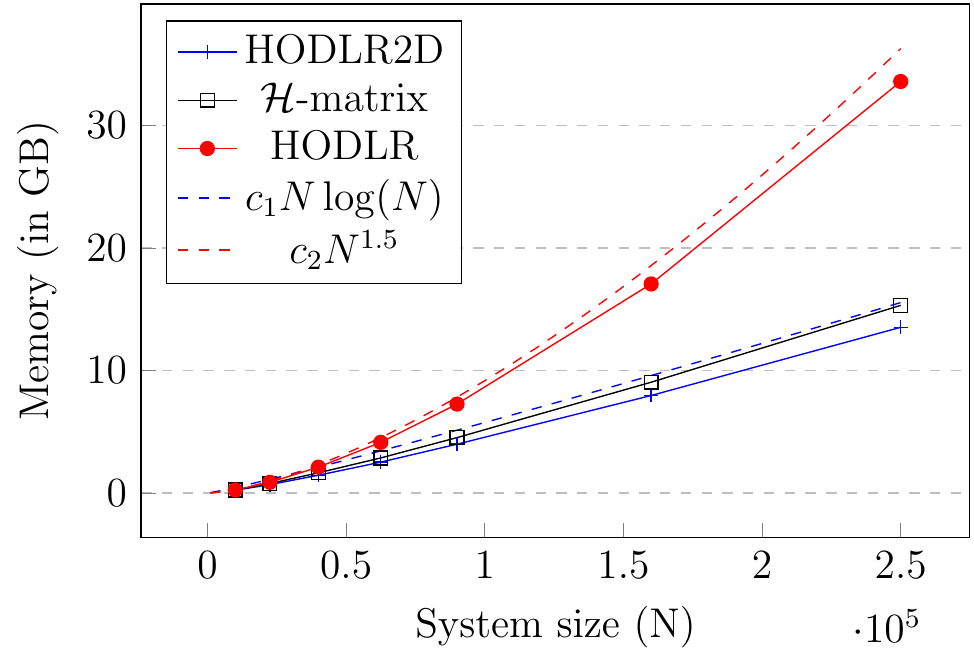}
		\caption{Storage}
		\label{fig:1rsub1}
	\end{subfigure}%
	\begin{subfigure}{.5\textwidth}
		\centering
		\includegraphics[width=.8\linewidth]{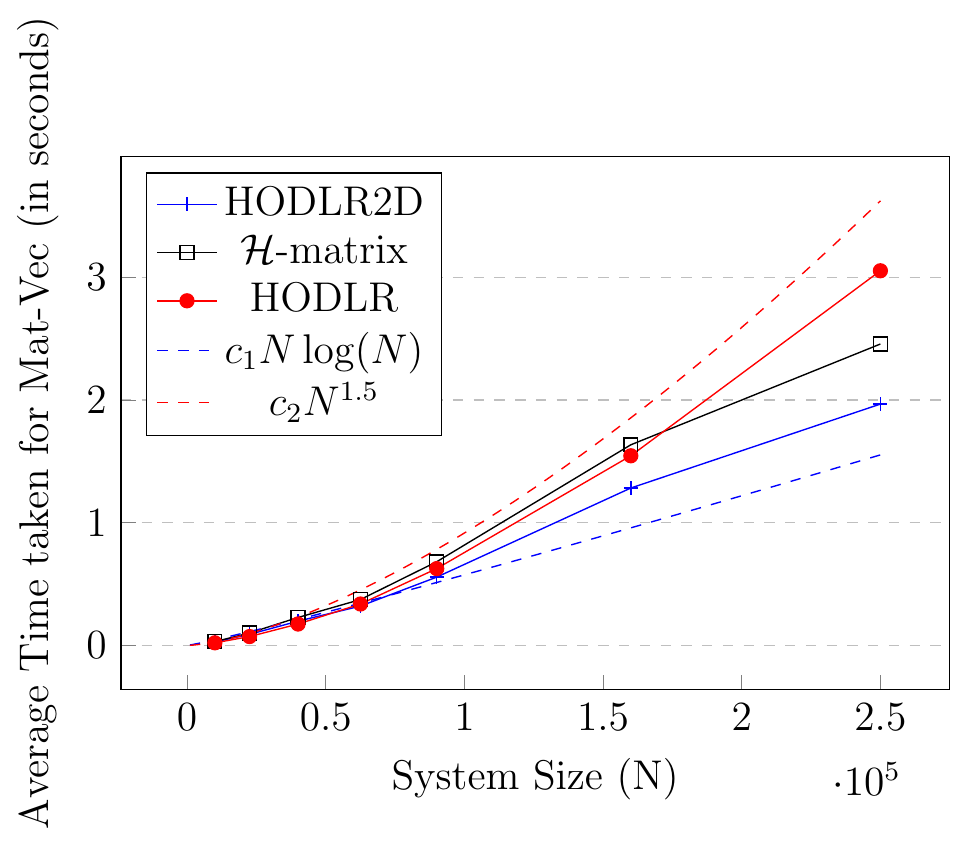}
		\caption{Average time taken for matrix-vector product}
		\label{fig:1rsub2}
	\end{subfigure}
	\caption{Performance comparison of HODLR2D, $\mathcal{H}$-matrix and HODLR for matrix $K$ whose elements are defined by Equation~\eqref{eq:1byr}}
	\label{fig:1r_perf}
\end{figure}

\subsection{HODLR2D accelerated iterative solver for radial basis function interpolation}

All the numerical experiments in this and the next subsection are performed without parallelization on a laptop with a $2.5$GHz Intel Core i5 processor and $16$GB RAM.

Let the location of the particles be a $\sqrt{N}\times \sqrt{N}$ Chebyshev grid in 2D over the domain $[-1,1]^2$. The dense linear system is generated using Equation~\eqref{eq:ker}.
\begin{equation} \label{eq:ker}
	\beta \lambda_{i} + \sum_{\substack{j=1 \\ j\neq i}}^{N}\phi(||x_i-x_j||_2)\lambda_{j} = f_i
\end{equation}
$\phi(r) : \mathbb{R}\rightarrow\mathbb{R}$, where $r$ is the Euclidean norm between two locations. The matrix form of Equation~\eqref{eq:ker} is given below:
\begin{equation} \label{eq:axb}
	A\lambda = f
\end{equation}

\begin{figure}[H]
	\centering
	\begin{subfigure}{.5\textwidth}
		\centering
		\includegraphics[width=0.8\linewidth]{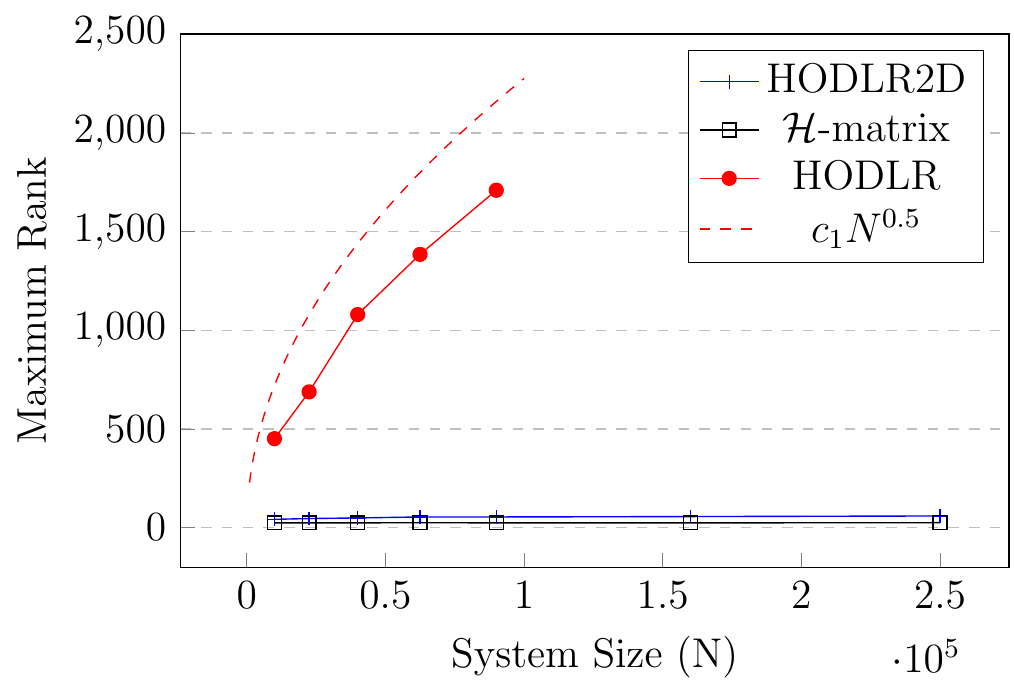}
		\caption{$r_m$ vs $N$ for $\phi_1(r)$}
		\label{fig:rank1}
	\end{subfigure}%
	\begin{subfigure}{.5\textwidth}
		\centering
		\includegraphics[width=0.8\linewidth]{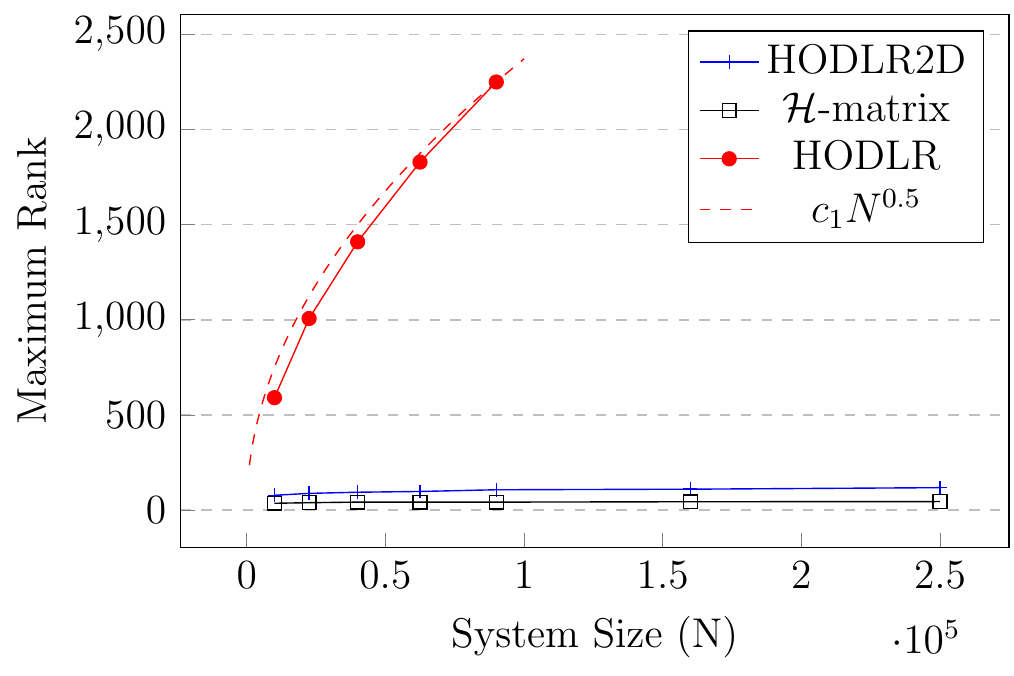}
		\caption{$r_m$ vs $N$ for $\phi_2(r)$}
		\label{fig:rank2}
	\end{subfigure}
	\caption{$r_m$ for HODLR2D, $\mathcal{H}$-matrix and HODLR for both radial basis functions}
	\label{fig:rank}
\end{figure}

\begin{figure}[H]
	\centering
	\begin{subfigure}{.5\textwidth}
		\centering
		\includegraphics[width=0.8\linewidth]{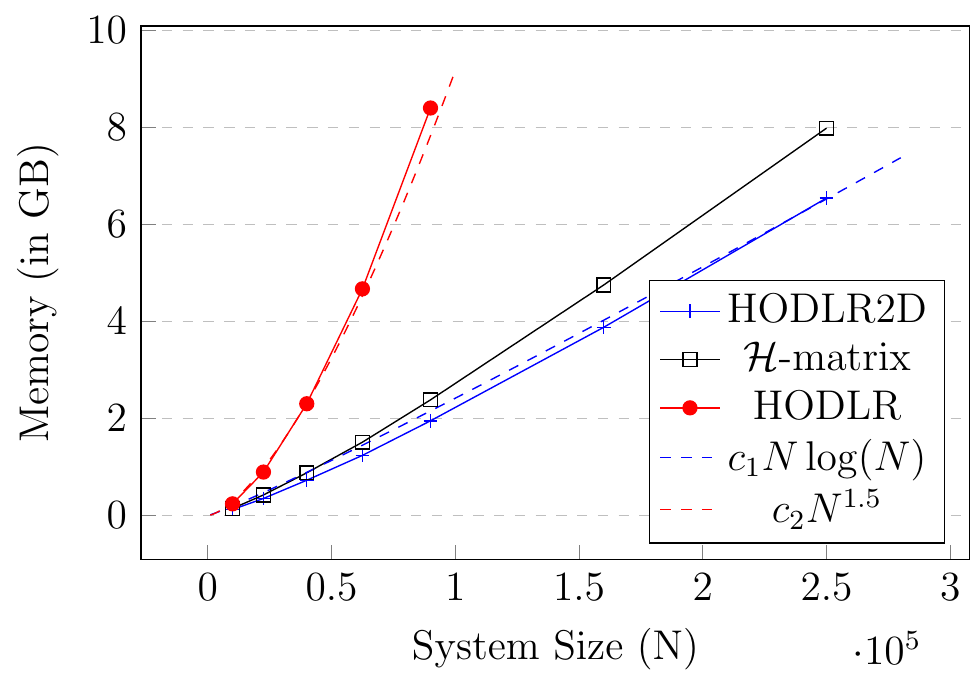}
		\caption{Storage in GB vs $N$ for $\phi_1(r)$}
		\label{fig:storage1}
	\end{subfigure}%
	\begin{subfigure}{.5\textwidth}
		\centering
		\includegraphics[width=0.8\linewidth]{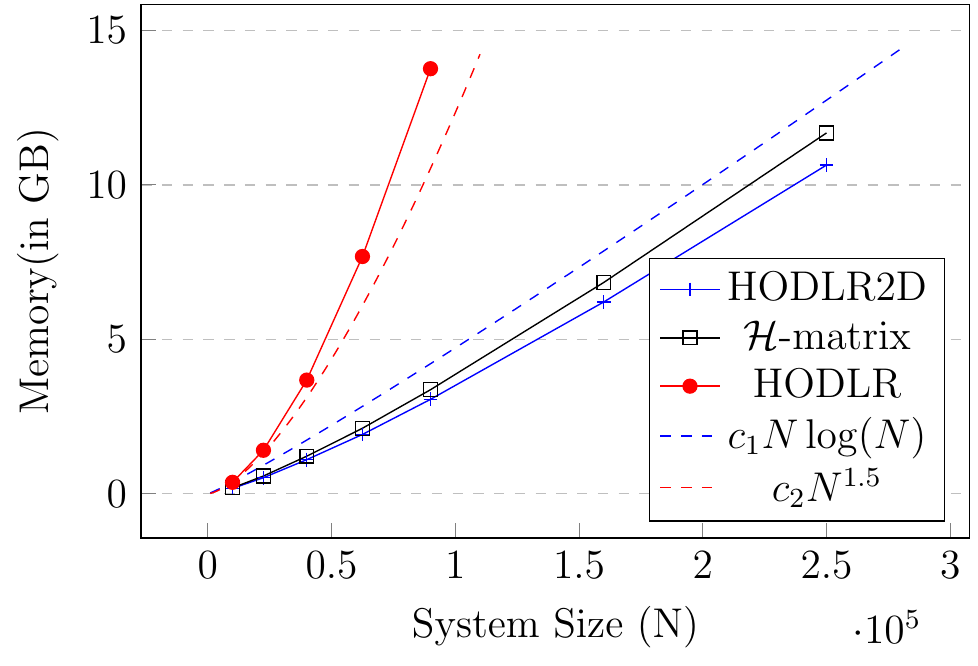}
		\caption{Storage in GB vs $N$ for $\phi_2(r)$}
		\label{fig:storage2}
	\end{subfigure}
	\caption{Storage for HODLR2D, $\mathcal{H}$-matrix and HODLR for both radial basis functions}
	\label{fig:storage}
\end{figure}

\begin{figure}[H]
	\centering
	\begin{subfigure}{.5\textwidth}
		\centering
		\includegraphics[width=0.8\linewidth]{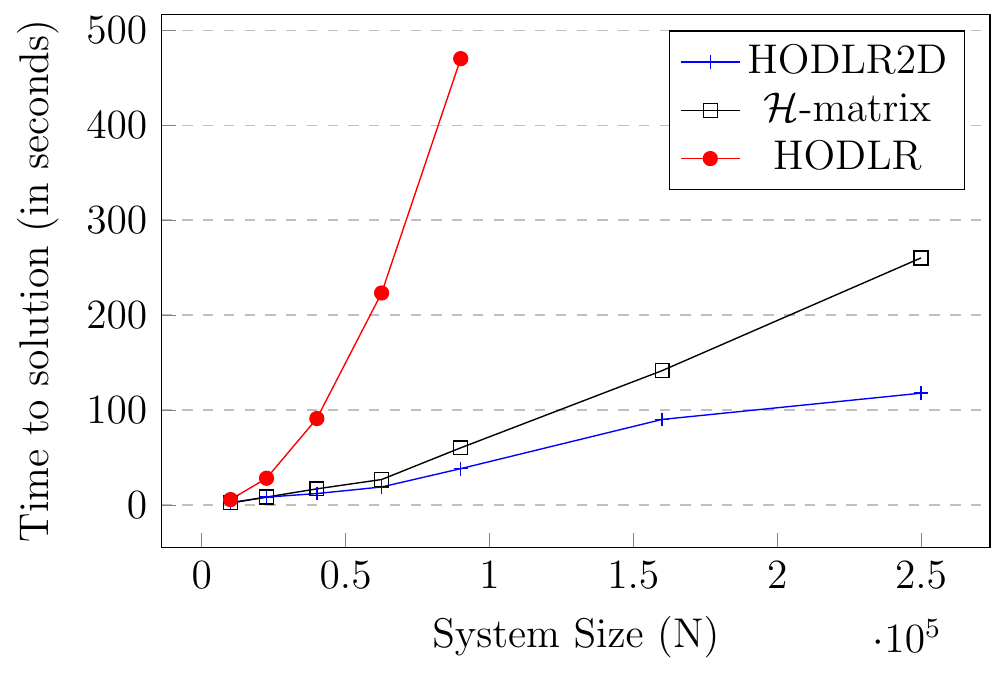}
		\caption{Total time in seconds vs $N$ for $\phi_1(r)$}
		\label{fig:time1}
	\end{subfigure}%
	\begin{subfigure}{.5\textwidth}
		\centering
		\includegraphics[width=0.8\linewidth]{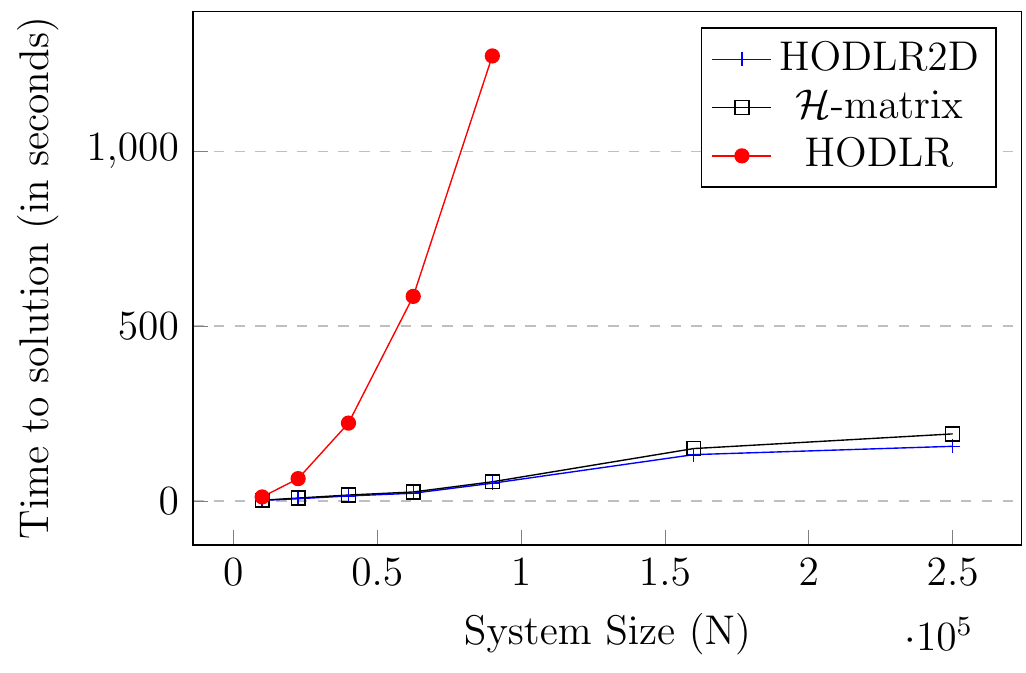}
		\caption{Total time in seconds vs $N$ for $\phi_2(r)$}
		\label{fig:time2}
	\end{subfigure}
	\caption{Time taken for the iterative solver using HODLR2D, $\mathcal{H}$-matrix and HODLR matrix-vector products for both radial basis functions}
	\label{fig:time}
\end{figure}

Two popularly used radial basis functions are considered.
\begin{equation}
	\text{Kernel $1$: }\phi_1(r) =
\begin{cases}
	\frac{\log_e{r}}{\log_e{a}} & \quad r\geq a\\
	\frac{r\log_e{r}-1}{a\log_e{a}-1} & \quad r<a
\end{cases}
\end{equation}
\begin{equation}
	\text{Kernel $2$: } \phi_2(r) =
		\begin{cases}
		\frac{a}{r} & \quad r\geq a\\
		\frac{r}{a} & \quad r<a
	\end{cases}
\end{equation}
With the choice of parameters $a=0.001$ and $\beta=N$, the resulting linear system corresponding to both the radial basis functions is well-conditioned. To verify the accuracy of the solution, a random vector $\lambda$ is used to produce the right-hand side vector $f$ in Equation~\eqref{eq:axb} by explicitly performing a matrix-vector product. We use the generated $f$ as our right-hand side for the iterative solver, and now we seek $\hat{\lambda}$. The relative error ($\epsilon_r$), measured in $2$-norm, i.e., $\dfrac{||\lambda-\hat{\lambda}||_2}{||\lambda||_2}$ is of the order of $10^{-10}$ in all the cases. The performance of different hierarchical formats with Kernel 1 is reported in Tables \ref{tab:ker1_init} and \ref{tab:k1_sol} and that for Kernel 2 in Tables \ref{tab:ker2_init} and \ref{tab:ker2_sol}. The scaling results are shown in Figures~\ref{fig:rank},~\ref{fig:storage} and \ref{fig:time}.
\begin{table}[H]
	\resizebox{\textwidth}{!}{
		\begin{tabular}{@{}cccccccccc@{}}\toprule
			\multirow{2}{*}{\textbf{N}} & \multicolumn{3}{c}{\textbf{$r_m$}}         & \multicolumn{3}{c}{\textbf{$T_I$}}  & \multicolumn{3}{c}{\textbf{CR}}    \\  \cmidrule(lr){2-4}\cmidrule(lr){5-7}\cmidrule(lr){8-10}
			& \textbf{HODLR2D} & \textbf{$\mathcal{H}$-matrix} & \textbf{HODLR} & \textbf{HODLR2D} & \textbf{$\mathcal{H}$-matrix} & \textbf{HODLR} & \textbf{HODLR2D} & \textbf{$\mathcal{H}$-matrix} & \textbf{HODLR} \\ \hline
			10000              & 42               & 24                            & 451            & 1.6              & 2                             & 4.8            & 0.154            & 0.182                         & 0.301          \\
			22500              & 46               & 24                            & 688            & 7.6              & 7.4                           & 25             & 0.087            & 0.105                         & 0.221          \\
			40000              & 49               & 24                            & 1080           & 10.7             & 15.2                          & 83.5           & 0.057            & 0.069                         & 0.18           \\
			62500              & 54               & 25                            & 1385           & 16.9             & 24.1                          & 206.5          & 0.04             & 0.048                         & 0.15           \\
			90000              & 54               & 24                            & 1710           & 33.8             & 50.1                          & 436.5          & 0.03             & 0.037                         & 0.129          \\
			160000             & 56               & 24                            & --             & 81.4             & 129                           & --             & 0.019            & 0.023                         & --             \\
			250000             & 59               & 25                            & --             & 105.4            & 231.8                         & --             & 0.013            & 0.016                         & --             \\ \hline
		\end{tabular}
	}
	\caption{Space complexity of HODLR2D, $\mathcal{H}$-matrix and HODLR for $\phi_1(r)$}
	\label{tab:ker1_init}
\end{table}

\begin{table}[H]
	\begin{tabular}{@{}cccc@{}}
		\toprule
		\multirow{2}{*}{$N$}  & \multicolumn{3}{c}{Matrix-Vector Product accelerated using}\\
		\cline{2-4} & \textbf{HODLR2D} & \textbf{$\mathcal{H}$-matrix} & \textbf{HODLR}\\
		\hline
		10000              & 0.231            & 0.288                         & 0.797 \\
		22500              & 0.556            & 0.852                         & 3.009 \\
		40000              & 1.198            & 1.716                         & 7.487 \\
		62500              & 1.961            & 2.559                         & 16.708 \\
		90000              & 4.32             & 10.028                        & 33.719 \\
		160000             & 8.643            & 12.451                        & --    \\
		250000             & 12.2752          & 28.44                         & --    \\ \hline
	\end{tabular}
	
	\caption{Time taken for the iterative solver $T_G$ to solve Equation~\eqref{eq:ker} for $\phi_1(r)$}
	\label{tab:k1_sol}
\end{table}

\begin{table}[H]
	\resizebox{\textwidth}{!}{
		\begin{tabular}{@{}cccccccccc@{}}\toprule
			\multirow{2}{*}{\textbf{N}} & \multicolumn{3}{c}{\textbf{$r_m$}}         & \multicolumn{3}{c}{\textbf{$T_I$}}  & \multicolumn{3}{c}{\textbf{CR}}    \\  \cmidrule(lr){2-4}\cmidrule(lr){5-7}\cmidrule(lr){8-10}
			& \textbf{HODLR2D} & \textbf{$\mathcal{H}$-matrix} & \textbf{HODLR} & \textbf{HODLR2D} & \textbf{$\mathcal{H}$-matrix} & \textbf{HODLR} & \textbf{HODLR2D} & \textbf{$\mathcal{H}$-matrix} & \textbf{HODLR} \\ \hline
			10000              & 78               & 36                            & 591            & 1.9              & 2.2                           & 11.3           & 0.218            & 0.239                         & 0.463          \\
			22500              & 88               & 39                            & 1007           & 6.4              & 8.2                           & 61.8           & 0.128            & 0.142                         & 0.347          \\
			40000              & 94               & 42                            & 1410           & 13.8             & 16.1                          & 214.8          & 0.086            & 0.095                         & 0.287          \\
			62500              & 98               & 42                            & 1829           & 21               & 24.2                          & 568.4          & 0.061            & 0.068                         & 0.246          \\
			90000              & 107              & 42                            & 2250           & 48.8             & 52.4                          & 1244.4         & 0.047            & 0.052                         & 0.213          \\
			160000             & 110              & 44                            & --             & 121.4            & 141.9                         & --             & 0.03             & 0.033                         & --             \\
			250000             & 118              & 45                            & --             & 147.6            & 178.9                         & --             & 0.021            & 0.023                         & --             \\ \hline
		\end{tabular}
	}
	\caption{Space complexity of HODLR2D, $\mathcal{H}$-matrix and HODLR for $\phi_2(r)$}
	\label{tab:ker2_init}
\end{table}

\begin{table}[H]
	\begin{tabular}{@{}cccc@{}}
		\toprule
		\multirow{2}{*}{$N$}  & \multicolumn{3}{c}{Matrix-Vector Product accelerated using}\\
		\cline{2-4} & \textbf{HODLR2D} & \textbf{$\mathcal{H}$-matrix} & \textbf{HODLR}\\
		\hline
		10000              & 0.142            & 0.168                         & 0.531 \\
		22500              & 0.36             & 0.515                         & 2.318 \\
		40000              & 0.896            & 1.129                         & 8.158 \\
		62500              & 1.465            & 1.799                         & 16.614 \\
		90000              & 2.221            & 2.769                         & 27.839 \\
		160000             & 11.336           & 8.195                         & --     \\
		250000             & 8.93127          & 13.06                         & --    \\ \hline
	\end{tabular}
	\caption{Time taken for the iterative solver $T_G$ to solve Equation~\eqref{eq:ker} for $\phi_2(r)$}
	\label{tab:ker2_sol}
\end{table}

\subsection{HODLR2D accelerated iterative solver for integral equations in \texorpdfstring{$2$}{2}D}
We now demonstrate the applicability of HODLR2D in solving the linear system arising out of discretization of Fredholm integral equation of the second kind.

Consider the integral equation in Equation~\eqref{eq:3}

\begin{equation} \label{eq:3}
	\psi(x) - \dfrac{i \kappa^{2}q(x)}4 \int_{\Omega} H_0^{(1)}\bkt{\kappa \magn{x-y}_2} \psi(y)dy = -\kappa^{2}q(x)\exp\bkt{i0.5x_1}
\end{equation}
where $x,y \in \Rb^2$, $x_1$ is the first coordinate of $x$, $\kappa=0.5$, $q(x)=1.5 \exp\bkt{-0.25 \magn{x}_2^2}$ and $\Omega = [-1,1]^2$.

Discretisation of Equation~\eqref{eq:3} is done as in~\cite{LS_2016,gujjula2022new} and this results in a linear system of the form
\begin{equation} \label{eq:5}
	A\vec{\psi} = \vec{f}
\end{equation}
where $\vec{\psi}$ is a vector of values of $\psi(x)$ at the grid points located in the leaf boxes of the quadtree as done in~\cite{LS_2016,gujjula2022new}.

\begin{table}[H]
	\resizebox{\textwidth}{!}{
		\begin{tabular}{@{}cccccccccc@{}}\toprule
			\multirow{2}{*}{\textbf{N}} & \multicolumn{3}{c}{\textbf{$r_m$}}         & \multicolumn{3}{c}{\textbf{$T_I$}}  & \multicolumn{3}{c}{\textbf{CR}}    \\  \cmidrule(lr){2-4}\cmidrule(lr){5-7}\cmidrule(lr){8-10}
			& \textbf{HODLR2D} & \textbf{$\mathcal{H}$-matrix} & \textbf{HODLR} & \textbf{HODLR2D} & \textbf{$\mathcal{H}$-matrix} & \textbf{HODLR} & \textbf{HODLR2D} & \textbf{$\mathcal{H}$-matrix} & \textbf{HODLR} \\ \hline
			1600               & 47               & --                            & 221            & 112.1            & 132.8                         & 101.7          & 0.807            & 1                             & 0.67           \\
			4096               & 42               & 17                            & 290            & 333.3            & 442.3                         & 420.1          & 0.344            & 0.452                         & 0.388          \\
			6400               & 51               & 17                            & 429            & 725.6            & 986.4                         & 973.3          & 0.318            & 0.429                         & 0.369          \\
			16384              & 44               & 16                            & 558            & 1967             & 2620.9                        & 3560.1         & 0.118            & 0.158                         & 0.198          \\
			25600              & 53               & 16                            & 830            & 4044.9           & 5564.6                        & 8269.4         & 0.104            & 0.143                         & 0.188          \\
			65536              & 47               & 15                            & 1067           & 10312            & 13808.7                       & 28797.6        & 0.037            & 0.049                         & 0.097          \\
			102400             & 55               & 15                            & --             & 20237.4          & 27950                         & --             & 0.031            & 0.043                         & --             \\ \hline
		\end{tabular}
	}
	\caption{Space Complexity of HODLR2D, $\mathcal{H}$-matrix and HODLR for Equation~\eqref{eq:5}}
	\label{tab:LS_init}
\end{table}

As done in~\cite{LS_2016,gujjula2022new}, we hierarchically subdivide the domain into smaller domains using a level restricted quadtree and represent the unknown function $\vec{\psi}(x)$ as polynomials on each of the leaf boxes of the level restricted quadtree. For our numerical experiments, we varied the number of levels in the level restricted quadtree and the number of points in the leaf level of the level restricted quadtree. \emph{It is to be noted that the level restricted quadtree mentioned here is in the context of discretizing the domain. The HODLR2D hierarchical low-rank structure still rests on the balanced quadtree as discussed in the previous section.}
\begin{table}[H]
	%	\resizebox{\textwidth}{!}{
	\begin{tabular}{@{}ccccccc@{}}\toprule
		\multirow{2}{*}{N} & \multicolumn{3}{c}{$T_G$ (in s)} & \multicolumn{3}{c}{Relative Error }                \\ \cmidrule(lr){2-4}\cmidrule(lr){5-7}
		& \textbf{HODLR2D} & \textbf{$\mathcal{H}$-matrix} & \textbf{HODLR} & \textbf{HODLR2D}      & \textbf{$\mathcal{H}$-matrix} & \textbf{HODLR}         \\ \hline
		1600               & 0.184            & 0.269                         & 0.421          & 0.567$\times10^{-10}$ & 0.613$\times10^{-10}$         & 0.616$\times10^{-10}$  \\
		4096               & 0.429            & 0.76                          & 1.91           & 0.634$\times10^{-10}$ & 0.608$\times10^{-10}$         & 0.62$\times10^{-10}$   \\
		6400               & 0.994            & 1.381                         & 3.649          & 0.613$\times10^{-10}$ & 0.636$\times10^{-10}$         & 0.632$\times10^{-10}$  \\
		16384              & 2.38             & 3.224                         & 14.678         & 0.758$\times10^{-10}$ & 0.604$\times10^{-10}$         & 0.674 $\times10^{-10}$ \\
		25600              & 5.217            & 7.146                         & 35.775         & 0.676$\times10^{-10}$ & 0.689$\times10^{-10}$         & 0.123$\times10^{-10}$  \\
		65536              & 11.612           & 15.617                        & 176.085        & 1.25$\times10^{-10}$  & 1.11$\times10^{-10}$          & 0.203$\times10^{-10}$  \\
		102400             & 24.7659          & 33.96                         & --             & 1.28$\times10^{-10}$  & 1.4 $\times10^{-10}$          & --                     \\ \hline
	\end{tabular}
	%		}
	\caption{Performance of HODLR2D, $\mathcal{H}$-matrix and HODLR for Equation~\eqref{eq:5}}
	\label{tab:ls_sol}
\end{table}

Since the right-hand side vector is well defined for this example, we use GMRES with stopping criteria as residual being $10^{-10}$ to get an approximate solution vector. The approximate solution vector is again used to generate the right-hand side by explicitly constructing the dense matrix and performing a matrix-vector product. The error reported in Table~\ref{tab:ls_sol} is the relative forward error in the solution vector $\psi$.

\begin{figure}[H]
	\centering
	\begin{subfigure}{.5\textwidth}
		\centering
		\includegraphics[width=0.8\linewidth]{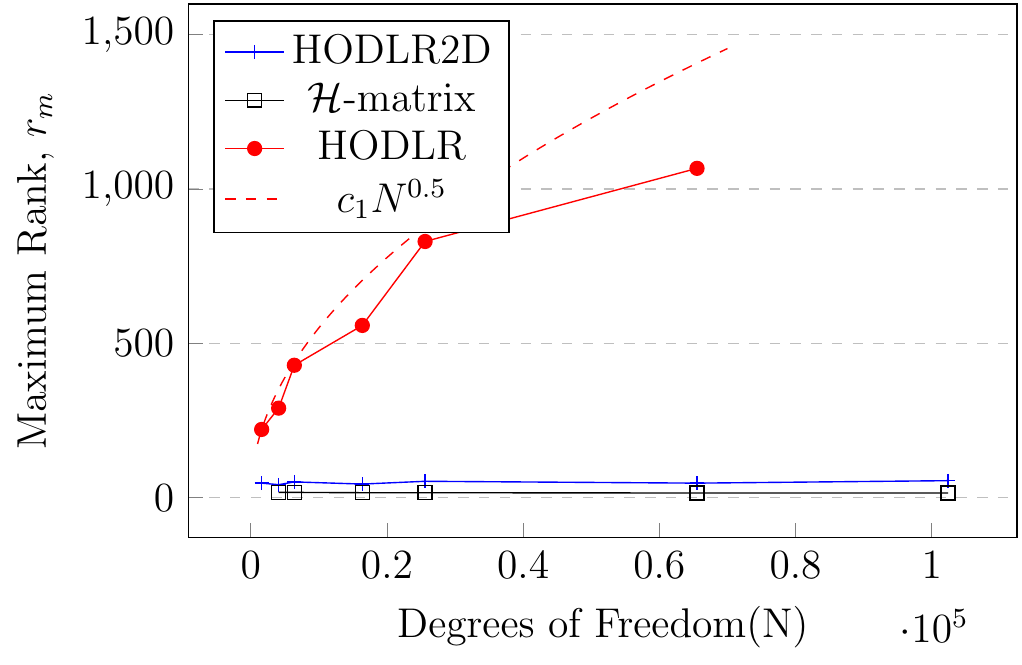}
		\caption{$r_m$ vs $N$}
		\label{fig:lssub1}
	\end{subfigure}%
	\begin{subfigure}{.5\textwidth}
		\centering
		\includegraphics[width=.8\linewidth]{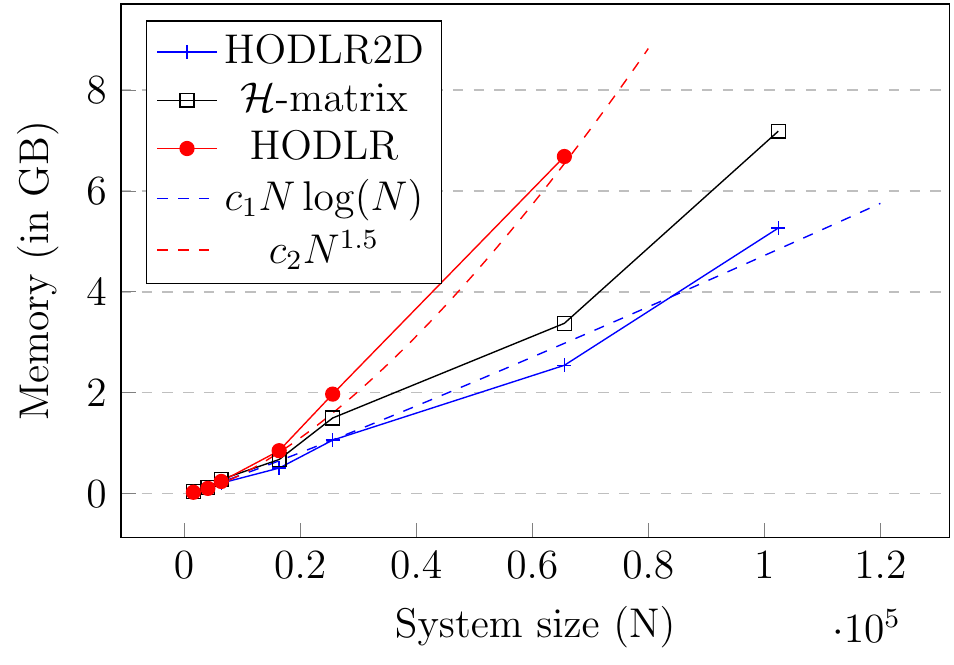}
		\caption{Memory(in GB) System Size(N)}
		\label{fig:lssub2}
	\end{subfigure}
	\caption{Storage for HODLR2D, $\mathcal{H}$-matrix and HODLR for Equation~\eqref{eq:5}}
	\label{fig:ls_mem}
\end{figure}

\subsection{Inferences}
The key observation from the Tables \ref{tab:ker1_init} through \ref{tab:ls_sol} is that the maximum rank reported remains almost constant in the case of  $\mathcal{H}$-matrix and HODLR2D, whereas for HODLR, it increases as a function of the system size. From Figures \ref{fig:1r_perf} through \ref{fig:ls_mem}, it is clear that HODLR2D scales almost linearly in terms of both storage and computational complexity. In all examples considered, HODLR2D performs significantly better in space and computational complexity than HODLR and provides an attractive alternative to $\mathcal{H}$-matrix.
	\section{Parallel HODLR2D}
\label{sec6}
This section will examine the parallel scalability of HODLR2D initialization (Algorithm~\ref{alg:cap}) and the matrix-vector product (Algorithm~\ref{alg:matvec}) on distributed memory systems. The existing literature on the parallel scalability of the rank structured matrices (both flat and hierarchical formats) is extensive, and we direct our readers to some seminal work\cite{parallel_Hmat,amestoy2019performance,akbudak2017tile,izadi2012hierarchical}. An important factor in the parallel scalability of an algorithm is how efficiently we can divide the load across the processes. In our case, we perform this by estimating the load of each node in the quadtree using Equation~\eqref{eq:par}.
\begin{equation}\label{eq:par}
	\text{Load} =
	\begin{cases}
		n^2 +( \abs{I_{C}} \times n) + \dsum_{i=1}^{\abs{I_{C}}}{m_i} +  n \times \bkt{ \dsum_{j=1}^{\abs{\mathcal{E}_{C}}}{k_j}},& \text{if leaf node}\\
		\abs{I_{C}} \times n + \dsum_{i=1}^{\abs{I_{C}}}{m_i}, & \text{if non-leaf node}
	\end{cases}
\end{equation}
where $n$ denotes the number of points in the cluster and $m_i$ denotes number of points in $i^{th}$ cluster in $I_C$ and finally $k_j$ denotes number of points in $j^{th}$ cluster in $\mathcal{E}_C$. Upon estimating the load for each node, we schedule the nodes across the processors by prioritizing the node with largest load first and by maintaining the load per processor almost constant. It is important to note that the scheduling is done apriori and not dynamically, the reason being it suits well for hybrid (multi-core + distributed) architectures. The initialization for HODLR2D does not require communication between the MPI processes. 

\begin{table}[H]
	\addtolength{\tabcolsep}{-1pt}
	\centering
	\begin{tabular}{@{}cccccccc@{}}
		\toprule
		\multirow{2}{*}{Number of MPI processes} & \multicolumn{7}{c}{System size(N)}                      \\ \cmidrule(l){2-8} 
		& 10000 & 22500 & 40000 & 62500 & 90000 & 160000 & 250000 \\ \midrule
		2 & 0.99  & 3.14  & 7.29  & 12.39 & 21.11 & 45.84  & 79.61  \\
		4 & 0.49  & 1.55  & 3.58  & 6.13  & 10.49 & 22.97  & 39.67  \\
		6 & 0.38  & 0.98  & 2.36  & 3.79  & 6.96  & 15.54  & 28.09  \\
		8 & 0.27  & 0.94  & 1.95  & 3.43  & 6.49  & 14.03  & 23.04  \\
		10 & 0.23  & 0.79  & 1.76  & 2.92  & 4.98  & 11.67  & 19.89  \\
		20 & 0.16  & 0.47  & 0.96  & 2.00  & 3.24  & 5.99   & 12.05  \\
		60 & 0.13  & 0.36  & 0.77  & 1.22  & 2.11  & 3.65   & 7.13   \\
		80 & 0.11  & 0.30  & 0.63  & 0.94  & 1.76  & 3.72   & 6.85  \\\bottomrule
	\end{tabular}
	\caption{HODLR2D Initialization time (in seconds) for different system sizes (N) by varying number of processors}
	\label{tab:par_ti}
\end{table}

\begin{table}[H]
	\addtolength{\tabcolsep}{-1pt}
	\centering
	\begin{tabular}{@{}cccccccc@{}}
		\toprule
		Number of MPI & \multicolumn{7}{c}{System size(N)}   \\ \cmidrule(l){2-8} 
		processes	& 10000 & 22500 & 40000 & 62500 & 90000 & 160000 & 250000 \\ \midrule
		2                                     & 0.016 & 0.048 & 0.106 & 0.172 & 0.298 & 0.671  & 1.034  \\
		4                                     & 0.008 & 0.025 & 0.056 & 0.091 & 0.154 & 0.346  & 0.530  \\
		6                                     & 0.007 & 0.021 & 0.043 & 0.076 & 0.110 & 0.253  & 0.383  \\
		8                                     & 0.005 & 0.015 & 0.030 & 0.049 & 0.087 & 0.197  & 0.290  \\
		10                                    & 0.005 & 0.014 & 0.028 & 0.046 & 0.073 & 0.159  & 0.243  \\
		20                                    & 0.004 & 0.009 & 0.021 & 0.031 & 0.048 & 0.105  & 0.162  \\
		60                                    & 0.016 & 0.016 & 0.015 & 0.015 & 0.032 & 0.050  & 0.086  \\
		80                                    & 0.012 & 0.011 & 0.015 & 0.019 & 0.017 & 0.032  & 0.051  \\ \bottomrule
	\end{tabular}
	\caption{HODLR2D Mat-Vec time (in seconds) for different system sizes (N) by varying number of processors}
	\label{tab:par_mv}
\end{table}
\begin{figure}[H]
	\centering
	\begin{subfigure}{.475\textwidth}
		\centering
		\includegraphics[width=.8\textwidth]{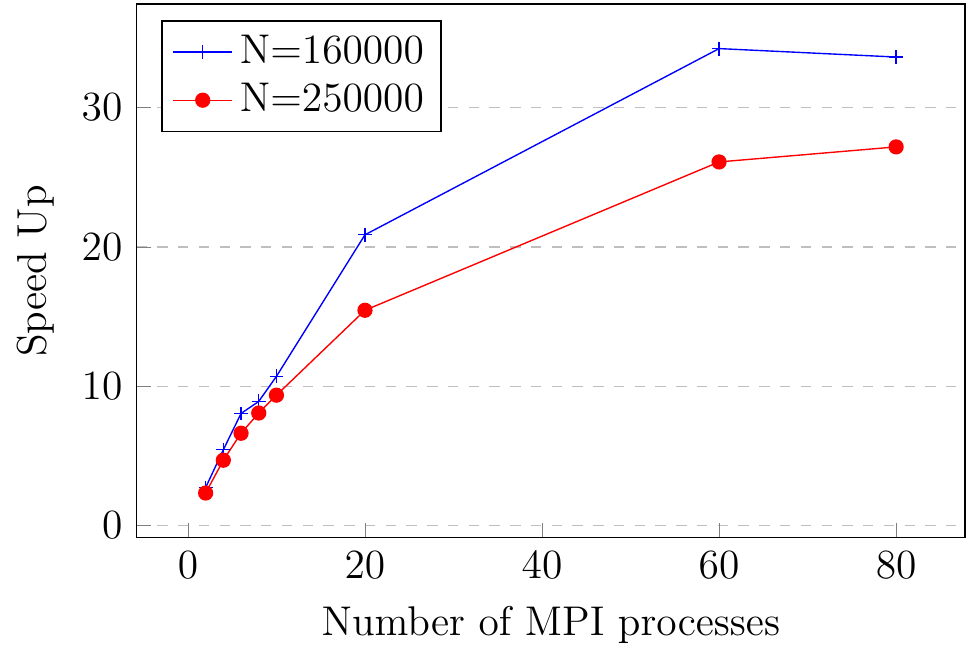}
		\caption{Speedup in Initialization time}
		\label{fig:partit}
	\end{subfigure}%
	\begin{subfigure}{.475\textwidth}
		\centering
		\includegraphics[width=.8\textwidth]{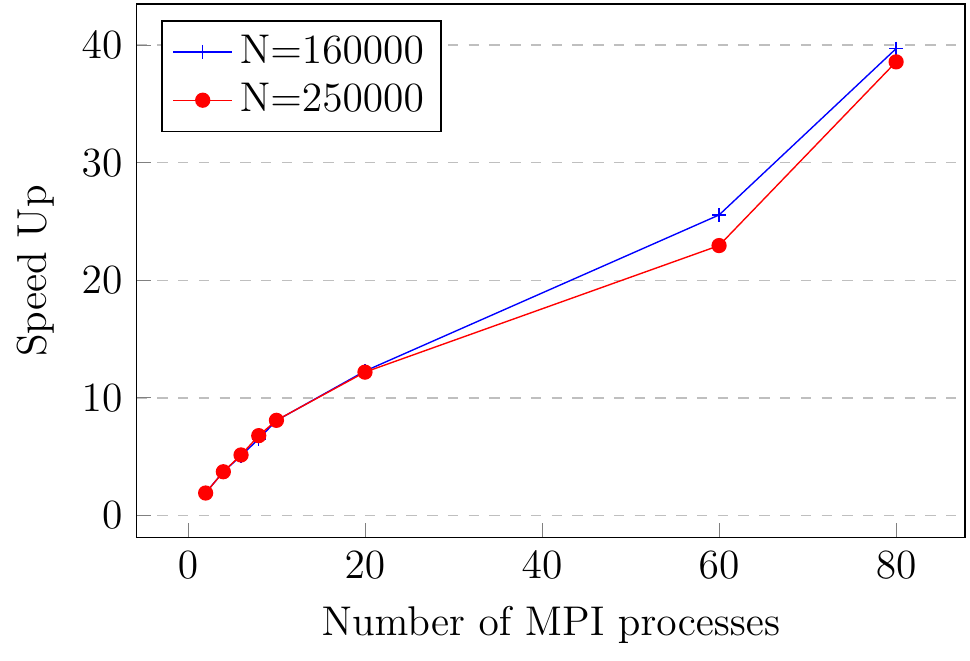}
		\caption{Speedup in Mat-Vec product time}
		\label{fig:par_mvt}
	\end{subfigure}
	\caption{Speedup of parallel HODLR2D vs number of MPI processors}
	\label{fig:par}
\end{figure}

Parallel HODLR2D is implemented using OpenMPI, and we repeat the same experiments in Section~\ref{sec_HODLR2D_matvec}, which serves as a baseline serial version of HODLR2D. We tabulate (Table~\ref{tab:par_ti} and Table~\ref{tab:par_mv}) the time to initialize and perform matrix-vector products using parallel HODLR2D by varying the number of MPI processes and system size. Figure \ref{fig:par} shows the speedup gained by parallel HODLR2D for particular system sizes.
	\section{Conclusion}
        We have presented a new hierarchical low-rank representation, HODLR2D, for a class of dense matrices arising out of two dimensional problems. We also provided theorems guaranteeing the growth of rank for different interactions in $2$D. These theorems form the basis of the HODLR2D algorithm. The key observation in HODLR2D is that the ranks of matrices corresponding to vertex sharing boxes are almost constant. We provide numerical benchmarks for our new HODLR2D structure by comparing the performance of HODLR2D with HODLR and $\mathcal{H}$-matrix with standard admissibility criterion. These benchmarks include time taken for matrix-vector products and to solve linear systems using GMRES based iterative solver. We also observe that both the memory requirement and time taken for matrix-vector products for HODLR2D matrix is significantly better than HODLR matrix. Further, the HODLR2D matrix structure provides an attractive alternative to $\mathcal{H}$-matrix. We also examine the parallel scalability of HODLR2D. We are also exploring on constructing a direct solver for this new HODLR2D structure and are also working on extending this to three dimensional problems.

\subsection*{Acknowledgments}
    We would like to acknowledge HPCE, IIT Madras for providing access to the AQUA cluster.
	
	% \bibliographystyle{unsrt}
	% \bibliography{references}

\end{document}